\documentclass[10pt,reqno]{amsart}

\usepackage[latin1]{inputenc}
\usepackage{amssymb,amsmath,amscd,amsfonts,color, amsthm}
\title[Asymptotic distribution of the Eigenvalues of a Gram Matrix]{The Empirical Distribution of the Eigenvalues of a Gram Matrix 
with a given variance profile}

\author[Hachem et al.]{W. Hachem, P. Loubaton and J. Najim}
\date{\today}

\newtheorem{theo}{Theorem}[section]

\newtheorem{lemma}[theo]{Lemma}
\newtheorem{coro}[theo]{Corollary}

\newtheorem{prop}[theo]{Proposition}
\newtheorem{assump}{Assumption A-\hspace{-0.15cm}}
\newcommand{\bdm}{\begin{displaymath}}
\newcommand{\edm}{\end{displaymath}}
\newcommand{\bea}{\begin{eqnarray*}}
\newcommand{\eea}{\end{eqnarray*}}

\newcommand{\Cplus}{\mathbb{C}^+}
\newcommand{\C}{\mathbb{C}}
\newcommand{\pitilde}{\tilde{\pi}}
\newcommand{\tv}{\mathrm{tv}}
\newcommand{\Cnabla}{\mathbb{C}^{\nabla}}
\newcommand{\im}{\mathrm{Im}}
\newcommand{\bs}{\boldsymbol}
\newcommand{\ti}{\tilde}
\newcommand{\Msub}{M_{\mathrm{sub}}}

\numberwithin{equation}{section}
\theoremstyle{remark}
\newtheorem{rem}{Remark}[section]

\begin{document}
\bibliographystyle{plain}
\begin{abstract}
  Consider a $N\times n$ random matrix $Y_n=(Y_{ij}^{n})$ where the
  entries are given by $Y_{ij}^{n}=\frac{\sigma(i/N,j/n)}{\sqrt{n}}
  X_{ij}^{n}$, the $X_{ij}^{n}$ being centered i.i.d.  and
  $\sigma:[0,1]^2 \rightarrow (0,\infty)$ being a continuous function
  called a variance profile.  Consider now a deterministic $N\times n$
  matrix $\Lambda_n=\left( \Lambda_{ij}^{n}\right)$ whose non diagonal
  elements are zero. Denote by $\Sigma_n$ the non-centered matrix $Y_n
  + \Lambda_n$ and by $N\wedge n= \min(N,n)$. Then under the
  assumption that $\lim_{n\rightarrow \infty } \frac Nn =c>0$ and
$$
\frac{1}{N\wedge n} \sum_{i=1}^{N\wedge n} 
\delta_{\left(\frac{i}{N\wedge n}, \left( \Lambda_{ii}^n\right)^2\right)} \xrightarrow[n\rightarrow \infty]{} H(dx,d\lambda),
$$
where $H$ is a probability measure, it is proven that
the empirical distribution of the eigenvalues of $
\Sigma_n \Sigma_n^T$ converges almost surely in distribution to a non random
probability measure. This measure is characterized in terms of its
Stieltjes transform, which is obtained with the help of an auxiliary system of equations.
This kind of results is of interest in the field of wireless communication.\\
\\
\noindent {\sc R\'esum\'e.} 
Soit $Y_n=(Y_{ij}^{n})$ une matrice $N\times n$ dont les entr\'ees sont donn\'ees par
$Y_{ij}^{n}=\frac{\sigma(i/N,j/n)}{\sqrt{n}} X_{ij}^{n}$, les $X_{ij}^{n}$ \'etant des variables al\'eatoires 
centr\'ees, i.i.d. et o\`u $\sigma:[0,1]^2 \rightarrow (0,\infty)$ est une fonction continue qu'on appellera profil de variance.
Consid\'erons une matrice d\'eterministe  $\Lambda_n=\left( \Lambda_{ij}^{n}\right)$ de dimensions $N\times n$ 
dont les \'el\'ements non diagonaux sont nuls. Appelons $\Sigma_n$ la matrice non centr\'ee d\'efinie par 
$\Sigma_n=Y_n + \Lambda_n$ et notons $N\wedge n=\min(N,n)$. Sous les hypoth\`eses que 
$\lim_{n\rightarrow \infty } \frac Nn =c>0$ et que 
$$
\frac{1}{N\wedge n} \sum_{i=1}^{N \wedge n}
\delta_{\left(\frac{i}{N\wedge n}, \left( \Lambda_{ii}^n\right)^2\right)} \xrightarrow[n\rightarrow \infty]{} H(dx,d\lambda),
$$
o\`u $H$ est une probabilit\'e, on d\'emontre que la mesure
empirique des valeurs propres de $\Sigma_n \Sigma_n^T$ converge
presque s\^urement vers une mesure de probabilit\'e d\'eterministe.
Cette mesure est caract\'eris\'ee par sa transform\'ee de Stieltjes,
qui s'obtient \`a l'aide d'un syst\`eme d'\'equations auxiliaire. Ce
type de r\'esultats pr\'esente un int\'er\^et dans le champ des
communications num\'eriques sans fil.
\end{abstract}

\maketitle

\section{Introduction} 
Consider a $N\times n$ random matrix $Y_n=(Y_{ij}^{n})$ where the entries are given by
\begin{equation}\label{variance-profile-variable}
Y_{ij}^{n}=\frac{\sigma(i/N,j/n)}{\sqrt{n}} X_{ij}^{n}
\end{equation}
where $\sigma: [0,1]\times [0,1] \rightarrow (0,\infty)$ is a
continuous function called a variance profile and the random variables
$X_{ij}^{n}$ are real, centered, independent and identically distributed
(i.i.d.) with finite $4+\epsilon$ moment. Consider a real deterministic
$N\times n$ matrix $\Lambda_n=(\Lambda^n_{ij})$ whose non-diagonal
elements are zero and consider the matrix $\Sigma_n=Y_n+\Lambda_n$. This model
has two interesting features: The random variables are independent but
not i.i.d. since the variance may vary and $\Lambda_n$, the centering perturbation of $Y_n$,
though (pseudo) diagonal can be of full rank. The purpose of this article is to study the convergence of the
empirical distribution of the eigenvalues of the Gram random matrix
$\Sigma_n \Sigma_n^T$ ($\Sigma_n^T$ being the transpose of $\Sigma_n$) when $n \rightarrow +\infty$ and $N \rightarrow +\infty$
in such a way that $\frac{N}{n} \rightarrow c$, $0 < c < +\infty$.  

The asymptotics of the spectrum of $N \times N$ Gram random matrices
$Z_n Z_n^{T}$ have been widely studied in the case where $Z_n$ is
centered (see Mar{\v{c}}enko and Pastur \cite{MarPas67}, Yin
\cite{Yin86}, Silverstein et al. \cite{Sil95,SilBai95}, Girko
\cite{Gir90-1,Gir01a}, Khorunzhy et al. \cite{KKP96}, Boutet de Monvel
et al.  \cite{BKV96}, etc.).  For an overview on asymptotic spectral
properties of random matrices, see Bai \cite{Bai99}. The case of a
Gram matrix $Z_n Z_n^T$ where $Z_n$ is non centered has comparatively
received less attention. Let us mention Girko (\cite{Gir01a}, chapter
7) where a general study is carried out for the matrix $Z_n = (W_n+
A_n)$ where $W_n$ has a given variance profile and $A_n$ is
deterministic. In \cite{Gir01a}, it is proved that the entries of the
resolvent $(Z_n Z_n^{T} - zI)^{-1}$ have the same asymptotic behavior
as the entries of a certain deterministic holomorphic $N \times N$
matrix valued function $T_n(z)$. This matrix-valued function is
characterized by a non linear system of $(n+N)$ coupled functional
equations. Using different methods, Brent Dozier and Silverstein
\cite{BreSil04pre} study the eigenvalue asymptotics of the matrix
$(R_n + X_n)(R_n+X_n)^T$ in the case where the matrices $X_n$ and
$R_n$ are independent random matrices, $X_n$ has i.i.d. entries and
the empirical distribution of $R_n R_n^T$ converges to a non-random
distribution. It is proved there that the eigenvalue distribution of
$(R_n + X_n)(R_n+X_n)^T$ converges almost surely towards a
deterministic distribution whose Stieltjes transform is uniquely
defined by a certain functional equation.  \\
\indent As in \cite{BreSil04pre}, the model studied in this article,
i.e.  $\Sigma_n= Y_n +\Lambda_n$, is a particular case of the general
case studied in (\cite{Gir01a}, chapter 7, equation $K_7$) for which
there exists a limiting distribution for the empirical distribution of
the eigenvalues. Since the centering term $\Lambda_n$ is
pseudo-diagonal, the proof of the convergence of the empirical
distribution of the eigenvalues is based on a direct analysis of the
diagonal terms of the resolvent $(\Sigma_n \Sigma_n^T -zI)^{-1}$. This
analysis leads in a natural way to the equations characterizing the
Stieltjes transform of the limiting probability distribution of the
eigenvalues.  
In the Wigner case with a variance profile (matrix $Y_n$
and the variance profile are symmetric), let us mention the recent
work of Anderson and Zeitouni \cite{AndZei05pre} where the asymptotics 
of the spectum and central limit theorems are investigated by means of systematic 
combinatorial enumeration.

Recently, many of these results have been applied to the field of
Signal Processing and Communication Systems and some new ones have
been developed for that purpose (Silverstein and Combettes
\cite{SilCom92}, Tse et al. \cite{TseHan99,TseZei00}, Debbah et al.
\cite{Debbah03}, Li et al. \cite{LTV04}, etc.).  The issue addressed
in this paper is mainly motivated by the performance analysis of
multiple-input multiple-output (MIMO) digital communication systems.
In MIMO systems with $n$ transmit antennas and $N$ receive antennas,
one can model the communication channel by a $N\times n$ matrix
$H_n=(H^n_{ij})$ where the entries $H^n_{ij}$ represent the complex
gain between transmit antenna $i$ and receive antenna $j$. The
statistics $ C_n = \frac{1}{n} \log \mbox{det} (I_n + \frac{H_n
  H_n^{\star}}{\sigma^{2}}) $ (where $H_n^{\star}$ is the hermitian adjoint and $\sigma^{2}$ represents the variance
of an additive noise corrupting the received signals) is a popular
performance analysis index since it has been shown in information
theory that $C_n$ is the maximum number of bits per channel use and per antenna
that can be transmitted reliably in a MIMO system with
channel matrix $H_n$. Since
$$
C_n= \frac{1}{n} \sum_{k=1}^{N} \log(1 + \frac{\mu_k}{\sigma^{2}}),
$$
where $(\mu_k)_{1\le k \le N}$ are the eigenvalues of $H_n
H_n^{\star}$, the empirical distribution of the eigenvalues of $H_n H_n^{\star}$
gives direct information on $C_n$ (see Tulino and Verdu
\cite{tul-ver-1} for an exhaustive review of recent results).  For
wireless systems, matrix $H_n$ is often modelled as a zero-mean Gaussian random
matrix and several articles have recently been
devoted to the study of the impact of the channel statistics (via the
eigenvalues of $H_n H_n^{\star}$) on the probability distribution of $C_n$
(Chuah et al. \cite{Chu02}, Goldsmith et al. \cite{Gol03}, see also
\cite{tul-ver-1} and the references therein).  Of particular interest
is also the channel matrix $H_n=F_N (Y_n +\Lambda_n) F_n^{T}$ where
$F_k=(F^k_{pq})_{1\le p,q\le k}$ is the Fourier matrix
(i.e. $F^k_{pq}=k^{-\frac 12} \exp(2 i \pi \frac{(p-1)(q-1)}{k})$) and
the matrix $Y_n$ is given by (\ref{variance-profile-variable}) (see
\cite{tul-ver-1}, p. 139 for more details). The matrices $H_n$ and $\Sigma_n$ having the same singular values, we 
will focus on the study of the empirical distribution of the singular values of $\Sigma_n$. Moreover, we will
focus on matrices with real entries since the complex case is a straightforward extension.

In the sequel, we will study simultaneously quantities (Stieltjes kernels) related to
the Stieltjes transforms of $\Sigma_n \Sigma_n^T$ and 
$\Sigma_n^T \Sigma_n$. 
Even if the Stieltjes transforms of $\Sigma_n \Sigma_n^T$ and $\Sigma_n^T \Sigma_n$
are related in an obvious way, the corresponding Stieltjes kernels are not, as we shall see.
We will prove that if
$N/n \xrightarrow[n\rightarrow \infty]{} c>0$ (since we study at the same time $\Sigma_n \Sigma_n^T$ and $\Sigma_n^T \Sigma_n$,
we assume without loss of generality that $c\le 1$) and if there
exists a probability measure $H$ on $[0,1]\times \mathbb{R}$ with compact support such
that
$$
\frac{1}{N} \sum_{i=1}^{N} 
\delta_{\left(\frac{i}{N}, \left(\Lambda_{ii}^n\right)^2\right)} \xrightarrow[n\rightarrow \infty]{{\mathcal D}} H(dx,d\lambda)
$$
where ${\mathcal D}$ stands for the convergence in distribution,
then almost surely, the empirical distribution of the eigenvalues of the random
matrix $\Sigma_n \Sigma_n^T$ (resp. $\Sigma^T_n \Sigma_n$) 
converges in distribution to a deterministic probability distribution $\mathbb{P}$
(resp. $\ti{\mathbb{P}}$). The probability distributions $\mathbb{P}$ and $\ti{\mathbb{P}}$ 
are characterized in terms of their Stieltjes transform
$$
f(z)=\int_{\mathbb{R}^+} \frac{ \mathbb{P}(dx)}{x-z}
\quad \textrm{and} \quad \ti f(z)=\int_{\mathbb{R}^+} \frac{\ti{ \mathbb{P}}(dx)}{x-z},\quad \textrm{Im}(z)\neq 0.
$$
as follows. Consider the following system of equations 
$$
\int g\,d\pi_z =\int \frac{g(u,\lambda)}{-z(1+\int \sigma^2(u,\cdot)d\pitilde_z) +
\frac {\lambda}{1+c \int \sigma^2(\cdot,c u) d\pi_z}} H(du,d\lambda)
$$
\begin{multline*}
\int g\,d\pitilde_z =c \int \frac{g(c u,\lambda)}{-z(1+c \int \sigma^2(\cdot,c u)d\pi_z) +
\frac {\lambda}{1+\int \sigma^2(u,\cdot) d\pitilde_z}} H(du,d\lambda)
\\ +(1-c) \int_c^1 \frac{g(u,0)}{-z(1+c \int \sigma^2(\cdot,u)d\pi_z)} \,du \\
\end{multline*}
where the unknown parameters are the complex measures $\pi_z$ and $\tilde{\pi}_z$. 
Both equalities stand for every continuous function $g:{\mathcal H} \rightarrow \mathbb{R}$, 
where ${\mathcal H}\subset [0,1]\times \mathbb{R}$ is the compact support of $H$. 
Then, this system admits a unique pair of solutions $(\pi_z(dx,d\lambda),\pitilde_z(dx,d\lambda))$.
In particular, $\pi_z$ is absolutely continuous with respect to $H$ while $\pitilde_z$ is not (see Remarks 
\ref{limiting-supports1}, \ref{limiting-supports2} and \ref{absolute-vodka} for more details). The Stieltjes transforms $f$ and $\ti f$  are then given by
$$
f(z)=\int_{[0,1]\times \mathbb{R}} \pi_z(dx,d\lambda)\quad \textrm{and}\quad \ti f(z)=\int_{[0,1]\times \mathbb{R}} 
\ti\pi_z(dx,d\lambda).
$$

The article is organized as follows. In Section \ref{resultats}, the notations and the assumptions are introduced and the main results
(Theorem \ref{existence-unicite} and Theorem \ref{convergence}) are stated. Section \ref{prooftheo1} is devoted to the proof 
of Theorem \ref{existence-unicite}. Section \ref{prooftheo2} is devoted to the proof 
of Theorem \ref{convergence}. In Section \ref{more-results}, we state as corollaries of Theorems \ref{existence-unicite}
and \ref{convergence} two simpler results which already exist in the literature. We also give details for the statement 
of the results in the complex setting (Section \ref{hypo-complex}).  

\section{Convergence of the Stieltjes Transform}\label{resultats}
\subsection{Notations and Assumptions}
Let $N=N(n)$ be a sequence of integers such that 
$$
\lim_{n\rightarrow \infty} \frac{N(n)}{n} =c.
$$
Consider a $N\times n$ random matrix $Y_n$ where the entries are given by 
$$
Y_{ij}^{n}=\frac{\sigma(i/N,j/n)}{\sqrt{n}} X_{ij}^{n}
$$
where $X_{ij}^n$ and $\sigma$ are defined below.
\begin{assump}\label{hypo1}
  The random variables $(X_{ij}^n\ ;\ 1\le i\le N,\,1\le j\le n\,,\, n\ge1)$ are real, independent and identically
  distributed. They are centered with $\mathbb{E}(X_{ij}^n)^2=1$ and
  satisfy:
$$
\exists\ \epsilon>0,\quad \mathbb{E}|X_{ij}^n|^{4+\epsilon}<\infty.
$$
\end{assump}
\noindent where $\mathbb{E}$ denotes the expectation.
\begin{assump}\label{hypo2}
  The real function $\sigma: [0,1] \times [0,1] \rightarrow \mathbb{R}$ is such that $\sigma^2$ is continuous.
  Therefore there exist a non-negative constant $\sigma_{\max}$ such that
\begin{equation}
\forall (x,y)\in [0,1]^2,\quad 0\le \sigma^2(x,y) \le \sigma^2_{\max}<\infty.
\end{equation}
\end{assump}

\begin{rem} As long as $\sigma^2$ remains continuous, $\sigma$ can vanish on portions of the domain $[0,1]^2$.
\end{rem}
Denote by $\mathrm{var}(Z)$ the variance of the random variable $Z$. Since $\mathrm{var}(Y_{ij}^n)=\sigma^2(i/N,j/n)/n$,
the function $\sigma$ will be called a variance profile.
Denote by $\delta_{z_0}(\,dz)$ the dirac measure at point $z_0$.
The indicator function of $A$ is denoted by $1_A(x)$.
Denote by $C_b({\mathcal X})$ (resp. $C_b({\mathcal X};\C))$ the set of real (resp. complex) 
continuous and bounded functions over the topological set ${\mathcal X}$ 
and by $\|f\|_{\infty}=\sup_{x\in {\mathcal X}} |f(x)|$, the supremum norm. If ${\mathcal X}$ 
is compact, we simply write $C({\mathcal X})$ (resp. $C({\mathcal X};\C))$) instead of  $C_b({\mathcal X})$
(resp. $C_b({\mathcal X};\C))$).

We will denote by $\xrightarrow[]{\mathcal D}$ the convergence in distribution for probability measures and by 
$\xrightarrow[]{w}$ the weak convergence for bounded complex measures.

Consider a real deterministic $N\times n$ matrix $\Lambda_n=(\Lambda_{ij}^n)$ whose non-diagonal entries are zero. 
We will often write $\Lambda_{ij}$ instead of $\Lambda_{ij}^n$.
We introduce the $N\times n$ matrix $\Sigma_n=Y_n + \Lambda_n$.  
For every matrix $A$, we will denote by $A^T$ its transpose, by $\mathrm{Tr}(A)$ its trace (if $A$ is square) 
and by $F^{A\,A^T}$, the empirical distribution function of 
the eigenvalues of $A\,A^T$. 
Since we will study at the same time the limiting spectrum of the matrices $\Sigma_n \Sigma_n^T$ and 
$\Sigma_n^T \Sigma_n$, we can assume without loss of generality that $c\le 1$. We also assume for simplicity that $N\le n$.

We assume that:
\begin{assump}\label{hypo3}
There exists a probability measure $H(\,du,d\lambda)$ over the set $[0,1]\times \mathbb{R}$ with compact support ${\mathcal H}$ such 
that
\begin{equation}
\frac1N \sum_{i=1}^N \delta_{\left( \frac iN,\,\left(\Lambda_{ii}^n\right)^2\right)}(du,d\lambda) \xrightarrow[n\rightarrow \infty]{\mathcal D} 
H(\,du,d\lambda).
\end{equation}
\end{assump}
\begin{rem}[The complex case] Assumptions (A-\ref{hypo1}), (A-\ref{hypo2}) and (A-\ref{hypo3}) must be slightly modified 
in the complex setting. Related modifications are stated in Section \ref{hypo-complex}.
\end{rem}

When dealing with vectors, the norm $\|\cdot\|$ will denote the Euclidean norm. In the case of matrices, the norm 
$\|\cdot\|$ will refer to the spectral norm. 

\begin{rem}\label{bornitude}
  Due to (A-\ref{hypo3}), we can assume without loss of generality
  that the $\Lambda_{ii}^n$'s are bounded for $n$ large enough. 
In fact, suppose not, then by
  (A-\ref{hypo3}), $\frac 1N \sum_{i=1}^N \delta_{\Lambda_{ii}^2}
  \rightarrow H_{\Lambda}(d\lambda)$ whose support is compact and,
  say, included in $[0,K]$. Then Portmanteau's theorem yields
  $\frac{1}{N} \sum_1^N 1_{[0,\ K+\delta]}(\Lambda_{ii}^2)
  \rightarrow 1$ thus
\begin{equation}\label{vanish}
\frac{\#\{ \Lambda_{ii},\ \Lambda_{ii}^2 \notin [0,\ K+\delta]\} }{N} 
= 1 - \frac{1}{N} \sum_1^N 1_{[0,\ K+\delta]}(\Lambda_{ii}^2) 
\xrightarrow[n\rightarrow \infty]{} 0.
\end{equation}
Denote by $\check{\Lambda}_n=( \check{\Lambda}_{ij}^n)$ the matrix whose non-diagonal elements are zero 
and set $\check{\Lambda}_{ii}^n 
= \Lambda_{ii}^n 1_{\{(\Lambda_{ii}^n)^2\le K+\delta\}}$. Then it is straightforward 
to check that 
$\frac1N \sum_{i=1}^N \delta_{\left( \frac iN,\,\check\Lambda_{ii}^2\right)}(du,d\lambda) \rightarrow
H(\,du,d\lambda)$. Moreover, if $\check \Sigma_n= Y_n+\check \Lambda_n$  then 
\begin{eqnarray*}
\| F^{\Sigma \Sigma^T} - F^{\check\Sigma \check\Sigma^T}\|_{\infty} \stackrel{(a)}{\le}
 \frac{\mathrm{rank}(\check\Sigma -\Sigma)}{N}
\stackrel{(b)}{\le}  \frac{\#\{ \Lambda_{ii},\ \Lambda_{ii}^2 \notin [\ 0,\ K+\delta]\} }{N} 
\xrightarrow[n\rightarrow \infty]{(c)} 0.
\end{eqnarray*} 
where (a) follows from Lemma 3.5 in \cite{Yin86} (see also
\cite{SilBai95}, Section 2), (b) follows from the fact that for a
rectangular matrix $A$, $\mathrm{rank(A)} \le$ the number of non zero
entries of $A$ and (c) follows from (\ref{vanish}). Therefore,
$F^{\Sigma \Sigma^T}$ converges iff $F^{\check\Sigma \check\Sigma^T}$
converges. In this case they share the same limit. Remark
\ref{bornitude} is proved.
\end{rem}

\begin{rem}\label{compact-support}
Due to Remark \ref{bornitude}, we will assume in the sequel that for all $n$, the support of 
$\frac1N \sum \delta_{\left( \frac iN,\,\Lambda_{ii}^2\right)}$ is included in a compact set 
${\mathcal K}\subset [0,1]\times \mathbb{R}$. 
\end{rem}

Let $\Cplus=\{z\in \C,\ \im(z)>0\}$ and
$\Cnabla=\{z\in \Cplus,\ |\mathrm{Re}(z)|\le \im(z)\}$.

\subsection{Stieltjes transforms and Stieltjes kernels}
Let $\nu$ be a bounded nonnegative measure over $\mathbb{R}$. Its Stieltjes transform $f$ is defined by:
$$
f(z)=\int_{\mathbb{R}} \frac{\nu(d\lambda)}{\lambda-z},\quad z\in\Cplus.
$$
We list below the main properties of the Stieltjes transforms that will be needed in the sequel.
\begin{prop} The following properties hold true:
\begin{enumerate}
\item Let $f$ be the Stieltjes transform of $\nu$, then
\begin{itemize}
\item[-] the function $f$ is analytic over $\Cplus$,
\item[-] the function $f$ satisfies: $|f(z)| \le \frac{\nu(\mathbb{R})}{\im(z)}$,
\item[-] if $z\in\Cplus$ then $f(z)\in \Cplus$,
\item[-] if $\nu(-\infty,0)=0$ then $z\in\Cplus$ implies $z\,f(z)\in \Cplus$.\\
\end{itemize}

\item Conversely, let $f$ be a function analytic over $\mathbb{C}^{+}$ such that $f(z) \in \mathbb{C}^{+}$ if 
$z \in \mathbb{C}^{+}$ and $|f(z)| |\mbox{Im}(z)|$ bounded on $\mathbb{C}^{+}$. Then, $f$ is the 
Stieltjes transform of a bounded positive measure $\mu$ and $\mu(\mathbb{R})$ is given by  
$$
\mu(\mathbb{R})= \lim_{y\rightarrow +\infty} -iy\,f(iy).
$$   
If moreover $z f(z) \in \mathbb{C}^{+}$ if $z \in \mathbb{C}^{+}$ then, $\mu(\mathbb{R}^-)=0$.\\

\item Let $\mathbb{P}_n$ and $\mathbb{P}$ be probability measures over $\mathbb{R}$ and denote by $f_n$ and $f$
their Stieltjes transforms. Then 
$$
\left( \forall z\in\Cplus,\ f_n(z) \xrightarrow[n\rightarrow\infty]{} f(z)\right) \quad \Rightarrow \quad 
\mathbb{P}_n\xrightarrow[n\rightarrow\infty]{\mathcal D} \mathbb{P}.  
$$

\end{enumerate}
\end{prop}

Let $A$ be an $n\times p$ matrix and let $I_n$ be the $n\times n$ identity. The resolvent of $A A^T$
is defined by 
$$
Q(z)=(A A^T -z\,I_n)^{-1}=\left(\rho_{ij}(z)\right)_{1\le i,j,\le n},\quad z\in \C\backslash\mathbb{R}.
$$
The following properties are straightforward.
\begin{prop}\label{proprietes-resolvante} Let $Q(z)$ be the resolvent of $A A^T$, then:
\begin{enumerate}
\item For all $z\in\Cplus$, $ \|Q(z)\|\le (\im(z))^{-1}$. Similarly,  $ |\rho_{ij}(z)|\le (\im(z))^{-1}$.\\
\item The function $h_n(z)=\frac1n \mathrm{Tr}\ Q(z)$ is the Stieltjes transform of the empirical 
distribution probability associated to the eigenvalues of $A A^T$. Since these eigenvalues are nonnegative, 
$z\,h_n(z)\in \Cplus$ for $z\in \Cplus$.\\
\item Let $\vec{\xi}$ be a $n\times 1$ vector, then $\im\left( z\vec{\xi} Q(z) \vec{\xi}^T\right) \in \Cplus$
for $z\in \Cplus$.
\end{enumerate}
\end{prop}

Denote by ${\mathcal M}_{\C}({\mathcal X})$ the set of complex measures over the topological set ${\mathcal X}$. 
In the sequel, we will call Stieltjes kernel every application 
$$
\mu: \Cplus  \rightarrow  {\mathcal M}_{\C}({\mathcal X})
$$
either denoted $\mu(z,dx)$ or $\mu_z(dx)$ and satisfying: 
\begin{enumerate}
\item\label{property1} $\forall g\in C_b({\mathcal X}),\ \int g\,d\mu_z$ is analytic over $\Cplus$,
\item\label{property2} $\forall z\in \Cplus,\ \forall g\in C_b({\mathcal X})$, 
$$
\left| \int g\,d\mu_z \right|\le \frac{\| g\|_{\infty}}{\im (z)}
$$
\item\label{property3} $\forall z\in \Cplus,\ \forall g\in C_b({\mathcal X})$ and $g\ge 0$ then 
$\im\left(\int g\,d\mu_z\right)\ge 0$,
\item\label{property4} $\forall z\in \Cplus,\ \forall g\in C_b({\mathcal X})$ and $g\ge 0$ then 
$\im\left( z\,\int g\,d\mu_z\right) \ge 0$.
\end{enumerate}

Let us introduce the following resolvents:
\begin{eqnarray*}
Q_n(z)&=&(\Sigma_n\Sigma_n^T -zI_N)^{-1}= \left( q_{ij}(z)\right)_{1\le i,j\le N},\quad z\in \Cplus,\\
\ti{Q}_n(z)&=& (\Sigma_n^T\Sigma_n -zI_n)^{-1}= \left(\ti{q}_{ij}(z)\right)_{1\le i,j\le n},\quad z\in \Cplus,
\end{eqnarray*}
and the following empirical measures defined for $z\in\Cplus$
\begin{eqnarray}
L^n_z(du,d\lambda) &=& \frac 1N \sum_{i=1}^N q_{ii}(z)\ \delta_{\left( \frac iN,\Lambda_{ii}^2\right)}(du,d\lambda),
\label{kernel1}\\
\ti{L}^n_z(du,d\lambda) &=& \frac 1n \sum_{i=1}^N \ti{q}_{ii}(z)\ \delta_{\left( \frac in, \Lambda_{ii}^2\right)}(du,d\lambda)
\nonumber\\
&\phantom{=}& + \left(\frac 1n \sum_{i=N+1}^n \ti{q}_{ii}(z)\ \delta_{\frac in}(du) \otimes \delta_0(d\lambda)\right)\, 1_{\{N<n\}},
\label{kernel2}
\end{eqnarray}
where $\otimes$ denotes the product of measures. Since $q_{ii}(z)$ (resp. $\ti q_{ii}(z)$) is analytic over $\Cplus$, satisfies 
$|q_{ii}(z)|\le (\im(z))^{-1}$ and $\min\left( \im(q_{ii}(z)), \im(zq_{ii}(z)) \right)>0$, $L^n$ (resp. $\ti L^n$) 
is a Stieltjes kernel. Recall that due to Remark \ref{compact-support}, $L^n$ and $\ti L^n$ have supports included in the compact 
set ${\mathcal K}$.

\begin{rem}[on the limiting support of $L^n$]\label{limiting-supports1}
Consider a converging subsequence of $L^n_z$,
then its limiting support is necessarily included in ${\mathcal H}$. 
\end{rem}

\begin{rem}[on the limiting support of $\ti L^n$]\label{limiting-supports2}
Denote by $H_c$ the image of the probability measure $H$ under the application $(u,\lambda)\mapsto (cu,\lambda)$, 
by ${\mathcal H}_c$ its support, by ${\mathcal R}$ the support of the measure 
$1_{[c,1]}(u)\,du \otimes \delta_0(d\lambda)$. Let $\ti {\mathcal H}={\mathcal H}_c \cup {\mathcal R}$.
Notice that $\ti {\mathcal H}$ is obviously compact.
Consider a converging subsequence of $\ti{L}^n_z$,
then its limiting support is necessarily included in $\ti {\mathcal H}$. 
\end{rem}

\subsection{Convergence of the empirical measures $L^n_z$ and $\ti{L}^n_z$}

\begin{theo}\label{existence-unicite}
Assume that $H$ is a probability measure over the set $[0,1]\times \mathbb{R}$ with compact support ${\mathcal H}$.
Assume moreover that (A-\ref{hypo2}) holds.
Then the system of equations
\begin{equation}\label{equation1}
\int g\,d\pi_z =\int \frac{g(u,\lambda)}{-z(1+\int \sigma^2(u,t)\pitilde(z,dt,d\zeta)) +
\frac {\lambda}{1+c \int \sigma^2(t,c u) \pi(z,dt,d\zeta)}} H(du,d\lambda)
\end{equation}
\begin{multline}\label{equation2}
\int g\,d\pitilde_z =c \int \frac{g(c u,\lambda)}{-z(1+c \int \sigma^2(t,c u)\pi(z, dt,d\zeta)) +
\frac {\lambda}{1+\int \sigma^2(u,t) \pitilde(z, dt,d\zeta)}} H(du,d\lambda)
\\ +(1-c) \int_c^1 \frac{g(u,0)}{-z(1+c \int \sigma^2(t,u)\pi(z,dt,d\zeta))} \,du \\
\end{multline}
where (\ref{equation1}) and (\ref{equation2}) hold for every $g\in C({\mathcal H})$, admits a unique couple of solutions
$(\pi(z,dt,d\lambda),\pitilde(z,dt,d\lambda))$ among the set of Stieltjes kernels for which the support of measure $\pi_z$ 
is included in ${\mathcal H}$ and the support of measure $\ti\pi_z$ is included in $\ti {\mathcal H}$.

Moreover the functions $f(z)=\int\,d\pi_z$ and $\ti f(z)=\int \,d\pitilde_z$ are the Stieltjes transforms 
of probability measures.
\end{theo}

\begin{rem}[on the absolute continuity of $\pi_z$ and $\pitilde_z$]\label{absolute-vodka}
  Due to (\ref{equation1}), the complex measure $\pi_z$ is absolutely
  continuous with respect to $H$.  However, it is clear
  from (\ref{equation2}) that $\pitilde_z$ has an absolutely
  continuous part with respect to $H_c$ (recall that $H_c$ is the
  image of $H$ under $(u,\lambda)\mapsto (cu,\lambda)$) and an
  absolutely continuous part with respect to $1_{[c,1]}(u)\,du \otimes
  \delta_0(d\lambda)$ (which is in general singular with respect to $H_c$).
  Therefore, it is much more convenient to work with Stieltjes kernels $\pi$ and $\pitilde$ 
  rather than with measure densities indexed by $z$.
\end{rem}

\begin{theo}\label{convergence}
Assume that (A-\ref{hypo1}), (A-\ref{hypo2}) and (A-\ref{hypo3}) hold and 
denote by $\pi$ and $\pitilde$ the two Stieltjes kernels solutions of the coupled equations (\ref{equation1}) and (\ref{equation2}).
Then 
\begin{enumerate}
\item Almost surely, the Stieltjes kernel $L^n$ defined by (\ref{kernel1}) 
converges weakly to $\pi$, that is:
  $$
  \textrm{a.s.}\quad \forall z\in \Cplus,
  \quad L^n_z
  \xrightarrow[n\rightarrow \infty]{w} \pi_z.
$$ 
\item Almost surely, the Stieltjes kernel $\ti{L}^n$ defined by (\ref{kernel2}) 
converges weakly to $\pitilde$.
\end{enumerate}
\end{theo}

Proofs of Theorems \ref{existence-unicite} and \ref{convergence} are  postponed to Sections
\ref{prooftheo1} and \ref{prooftheo2}.

\begin{coro}\label{convergence-distrib}Assume that (A-\ref{hypo1}), (A-\ref{hypo2}) and (A-\ref{hypo3}) hold and denote by $\pi$ 
and $\pitilde$ the two Stieltjes kernels solutions of the coupled equations (\ref{equation1}) and (\ref{equation2}).
Then the empirical distribution of the eigenvalues of the matrix $\Sigma_n \Sigma_n^T$ converges almost surely to 
a non-random probability measure $\mathbb{P}$ whose Stieltjes transform $f(z)=\int_{\mathbb{R}^+} \frac{ \mathbb{P}(dx)}{x-z}$
is given by:
$$
f(z)=\int_{\mathcal H} \pi_z(dx,d\lambda).
$$
Similarly, the empirical distribution of the eigenvalues of the matrix $\Sigma_n^T  \Sigma_n$ converges almost surely to 
a non-random probability measure $\ti{\mathbb{P}}$ whose Stieltjes transform $\ti f(z)$ 
is given by:
$$
\ti{f}(z)=\int_{\ti{\mathcal H}} \ti \pi_z(dx,d\lambda).
$$
\end{coro}

\begin{proof}[Proof of Corollary \ref{convergence-distrib}]
The Stieltjes transform of $\Sigma_n \Sigma_n^T$ is equal to $\frac 1N \sum_{i=1}^N q_{ii}(z) =\int  \,dL_z^n$. By Theorem 
\ref{convergence},
\begin{equation}\label{conv-transfo-stieljes}
\mathrm{a.s.}\quad \forall z\in \Cplus,\quad \int \,dL_z^n \xrightarrow[n\rightarrow\infty]{} \int \,d\pi_z.
\end{equation}
Since $\int\,d \pi_z$ is the Stieltjes transform of a probability measure $\mathbb{P}$ by Theorem \ref{existence-unicite},
eq. (\ref{conv-transfo-stieljes}) implies that $F^{\Sigma_n \Sigma_n^T}$ converges weakly to $\mathbb{P}$.
One can similarly prove that $F^{\Sigma_n^T \Sigma_n}$ converges weakly to a probability measure $\mathbb{\ti P}$.
\end{proof}

\section{Proof of Theorem \ref{existence-unicite}}\label{prooftheo1}
We first introduce some notations.
Denote by 
\begin{eqnarray*}
D(\pitilde_z,\pi_z)(u,\lambda)& =& -z\left(1+\int \sigma^2(u,t)\pitilde(z,dt,d\zeta)\right) +
\frac {\lambda}{1+c \int \sigma^2(t,c u) \pi(z,dt,d\zeta)},\\
d(\pi_z)(u) &=& 1+c \int \sigma^2(t,c u) \pi(z,dt,d\zeta),\\
\tilde{D}(\pi_z,\pitilde_z)(u,\lambda) &=&-z\left(1+c \int \sigma^2(t,c u)\pi(z, dt,d\zeta)\right) +
\frac {\lambda}{1+\int \sigma^2(u,t) \pitilde(z, dt,d\zeta)},\\
\tilde{d}(\pitilde_z)(u)&=& 1+\int \sigma^2(u,t) \pitilde(z, dt,d\zeta),\\
\kappa(\pi_z)(u)&=&-z\left(1+c \int \sigma^2(t,u)\pi(z,dt,d\zeta)\right).
\end{eqnarray*}
Let $\nu$ be a complex measure over the set ${\mathcal H}$ (recall that ${\mathcal H}$ is compact by
(A-\ref{hypo3}))
then we denote by $\|\nu\|_{\tv}$ the total variation 
norm of $\nu$, that is 
\begin{eqnarray*}
\|\nu\|_{\tv}& = & |\nu|({\mathcal H})\\
&=& \sup \left\{\left|\int f \,d\nu\right|,\ f\in C({\mathcal H};\C),\ \|f\|_{\infty} \le 1\right\}.
\end{eqnarray*}
\subsection{Proof of the unicity of the solutions}
Notice that the system of equations (\ref{equation1}) and (\ref{equation2}) remains
true for every $g\in C({\mathcal H};\C)$ (consider $g=h+ik$) and assume that both 
$(\pi_z,\pitilde_z)$ and $(\rho_z,\tilde{\rho}_z)$ are pairs of solutions of the given system.  
Let $g\in C({\mathcal H};\C)$, 
then (\ref{equation1}) yields:
\begin{eqnarray*}
\int g\,d \pi_z -\int g\,d\rho_z
& = & \int \frac{z\, g(u,\lambda) \int \sigma^2(u,t)(\pitilde(z,dt,d\zeta)-\tilde{\rho}(z,dt,d\zeta))}
{D(\pitilde_z,\pi_z)\times D(\tilde{\rho}_z,\rho_z)}H(du,d\lambda)\\
&\phantom{=}& + \int \frac{c \lambda\,g(u,\lambda)\int \sigma^2(t,cu) (\rho(z,dt,d\zeta)-\pi(z,dt,d\zeta))}
{D(\pitilde_z,\pi_z)\times D(\tilde{\rho}_z,\rho_z)\times d(\pi_z)\times  d(\rho_z)} H(du,d\lambda)\\
\end{eqnarray*}
and
\begin{eqnarray*}
\left| \int g\,d \pi_z -\int g\,d\rho_z \right| 
&\le& |z| \sigma^2_{\max}\,\|g\|_{\infty}\, \|\pitilde_z -\tilde{\rho}_z\|_{\tv} \int\frac{dH}
{\left|D(\pitilde_z,\pi_z)\times D(\tilde{\rho}_z,\rho_z)\right|}\\
&\phantom{=}& +\,  c\, \sigma^2_{\max}\,\|g\|_{\infty}\, \|\pi_z -\rho_z\|_{\tv} \int \frac{\lambda H(du,d\lambda)}
{\left| D(\pitilde_z,\pi_z)\times D(\tilde{\rho}_z,\rho_z)\times d(\pi_z)\times  d(\rho_z)\right|} \\
\end{eqnarray*}
If one takes the supremum over the functions $g\in C({\mathcal H};\C), \|g\|_\infty\le 1$, one gets :
$$
\| \pi_z -\rho_z \|_{\tv} \le {\boldsymbol \alpha}\|\pi_z -\rho_z\|_{\tv} 
+ {\boldsymbol \beta} \|\pitilde_z -\tilde{\rho}_z\|_{\tv}
$$
where 
\begin{eqnarray*}
{\boldsymbol \alpha}={\boldsymbol \alpha}(\pi,\rho,\pitilde,\tilde{\rho})&=& c\, \sigma^2_{\max} \int \frac{\lambda H(du,d\lambda)}
{\left| D(\pitilde_z,\pi_z)\times D(\tilde{\rho}_z,\rho_z)\times d(\pi_z)\times  d(\rho_z)\right|} \\
{\boldsymbol \beta}={\boldsymbol \beta}(\pi,\rho,\pitilde,\tilde{\rho})&=& |z| \sigma^2_{\max}\int\frac{dH}
{\left| D(\pitilde_z,\pi_z)\times D(\tilde{\rho}_z,\rho_z)\right|}\\
\end{eqnarray*}
Similarly, (\ref{equation2}) yields:
\begin{eqnarray*}
\int g\,d \pitilde_z -\int g\,d\tilde{\rho}_z
& = & c \int \frac{c\,z\, g(u,\lambda) \int \sigma^2(t,cu)(\pi(z,dt,d\zeta)-\rho(z,dt,d\zeta))}
{\tilde{D}(\pi_z,\pitilde_z)\times \tilde{D}(\rho_z,\tilde{\rho}_z)}H(du,d\lambda)\\
&\phantom{=}& + c \int \frac{\lambda\,g(u,\lambda)\int \sigma^2(u,t) (\pitilde(z,dt,d\zeta)-\tilde{\rho}(z,dt,d\zeta))}
{\tilde{D}(\pi_z,\pitilde_z)\times \tilde{D}(\rho_z,\tilde{\rho}_z)\times \tilde{d}({\pitilde}_z)\times \tilde{d}(\tilde{\rho}_z)}
H(du,d\lambda)\\
&\phantom{=}& \phantom{+} +(1-c) \int_c^1\frac{c\,z\,g(u,0)\int \sigma^2(t,u)(\pi(z,dt,d\zeta)-\rho(z,dt,d\zeta))}
{\kappa(\pi_z)\times \kappa(\rho_z)}du
\end{eqnarray*}
and
\begin{eqnarray*}
\left| \int g\,d \pitilde_z -\int g\,d\tilde{\rho}_z\right| 
&\le & c^2 \sigma_{\max}^2 |z|\,\| g\|_{\infty} \|\pi_z -\rho_z \|_{\tv} \int \frac{dH}
{\left| \tilde{D}(\pi_z,\pitilde_z)\times \tilde{D}(\rho_z,\tilde{\rho}_z)\right|}\\
&\phantom{=}& +\ c\,\sigma_{\max}^2  \| g\|_{\infty} \|\pitilde_z -\tilde{\rho}_z \|_{\tv}
\int \frac{\lambda H(du,d\lambda)}
{\left|\tilde{D}(\pi_z,\pitilde_z)\times \tilde{D}(\rho_z,\tilde{\rho}_z)\times \tilde{d}(\pitilde_z)\times \tilde{d}(\tilde{\rho}_z)\right|}\\
&\phantom{=}& \phantom{+} +(1-c)c \sigma_{\max}^2 |z|\,\| g\|_{\infty} \|\pi_z -\rho_z \|_{\tv}
\int_c^1 \frac{du}{\left|\kappa(\pi_z)\times \kappa(\rho_z)\right|}
\end{eqnarray*}
As previously, by taking the supremum over $g\in C({\mathcal H};\C), \|g\|_\infty\le 1$, we get:
$$
\|\pitilde_z -\tilde{\rho}_z\|_{\tv} \le \boldsymbol{\tilde{\alpha}}\|\pi_z -\rho_z\|_{\tv} 
+ \boldsymbol{\tilde{\beta}} \|\pitilde_z -\tilde{\rho}_z\|_{\tv}
$$
where 
\begin{eqnarray*}
\boldsymbol{\tilde{\alpha}}= \boldsymbol{\tilde{\alpha}}(\pi,\rho,\pitilde,\tilde{\rho})
&=& c\, \sigma^2_{\max} |z| \left( c \int \frac{dH}
{\left| \tilde{D}(\pi_z,\pitilde_z)\times \tilde{D}(\rho_z,\tilde{\rho}_z)\right|}
+ (1-c) \int_c^1 \frac{du}{\left| \kappa(\pi_z)\times \kappa(\rho_z)\right|} \right)\\
\boldsymbol{\tilde{\beta}}= \boldsymbol{\tilde{\beta}}(\pi,\rho,\pitilde,\tilde{\rho})&=& 
c\, \sigma^2_{\max}
\int \frac{\lambda H(du,d\lambda)}
{\left|\tilde{D}(\pi_z,\pitilde_z)\times \tilde{D}(\rho_z,\tilde{\rho}_z)\times \tilde{d}(\pitilde_z)\times \tilde{d}(\tilde{\rho}_z)
\right|}\\
\end{eqnarray*}
We end up with the following inequations:
\begin{equation}\label{system}
\left\{
\begin{array}{ll}
\| \pi_z -\rho_z \|_{\tv} & \le {\boldsymbol \alpha}\|\pi_z -\rho_z\|_{\tv} 
+ {\boldsymbol \beta} \|\pitilde_z -\tilde{\rho}_z\|_{\tv}\\
\|\pitilde_z -\tilde{\rho}_z\|_{\tv} & \le \boldsymbol{\tilde{\alpha}}\|\pi_z -\rho_z\|_{\tv} 
+ \boldsymbol{\tilde{\beta}} \|\pitilde_z -\tilde{\rho}_z\|_{\tv}
\end{array}
\right.
\end{equation}
Let us prove now that for $z\in \Cnabla$ with $\im(z)$ large enough, then $\boldsymbol{\alpha}<\frac 12$. 

Since $\pi$ and $\pitilde$ are assumed to be Stieltjes kernels, $\im\left(z\int \sigma^2(u,t)\pitilde(z,dt,d\zeta)\right)\ge 0$
and $\im\left(\int \sigma^2(t,c u) \pi(z,dt,d\zeta)\right)\ge 0$. Therefore, $ \im(D(\pitilde_z,\pi_z)) \leq -\im(z)$ and hence
$|\im(D(\pitilde_z,\pi_z))|\ge \im(z)$. Similarly, 
$|\im(D(\tilde{\rho}_z,\rho_z))|\ge \im(z)$. Thus,
$$
\frac{1}{|D(\pitilde_z,\pi_z)\times D(\tilde{\rho}_z,\rho_z)|}\le \frac{1}{\im^2(z)}.
$$
Now consider $z d(\pi_z)$. As previously, $\im(z d(\pi_z)) \ge \im(z)$.
As $|z d(\pi_z)| \geq |\im(z d(\pi_z))|$, this implies that $\frac{1}{|z d(\pi_z)|}\le \frac{1}{\im(z)}$ and
$\frac{1}{|d(\pi_z)|}\le \frac{|z|}{\im(z)}$. Since $z \in \Cnabla$,
$\frac{|z|}{\im(z)}\le \sqrt{2}$. The same argument holds for 
$d(\rho_z)$ thus we get
$$
{\boldsymbol \alpha} \le \frac{2\,c\, \sigma^2_{\max} \int \lambda H(du,d\lambda)}{\im^2(z)}<\frac 12 \quad 
\mathrm{for}\  \im(z)\ \mathrm{large\ enough.}
$$
With similar arguments, one can prove that
\begin{equation}\label{majoration}
\bs{\beta} \le \frac{\sqrt{2}\sigma^2_{\max}}{\im(z)},\quad \bs{\ti{\alpha}} \le \frac{3 \sigma^2_{\max}}{\im(z)},
\quad \bs{\ti{\beta}} \le \frac{2 \sigma_{\max}^2}{\im^2(z)} \int \lambda H(\,dud\lambda).
\end{equation}

Therefore $\max(\bs{\alpha},\boldsymbol{\beta}, \boldsymbol{\tilde{\alpha}},\boldsymbol{\tilde{\beta}})\le \theta<\frac 12$
for $z\in \Cnabla$ and $\im(z)$ large enough where $\theta$ does not depend on $(\pi, \ti{\pi}, \rho, \ti{\rho})$. 
Therefore, the system (\ref{system}) yields 
$$
\|\pi_z -\rho_z \|_{\tv}= \|\pitilde_z
-\tilde{\rho}_z\|_{\tv}=0\quad \mathrm{for}\ z\in \Cnabla\
\textrm{and}\ \im(z)\ \mathrm{large\ enough.}
$$
Now take $z\in \Cplus$ and $g \in C({\mathcal H})$.
Since $\int g \,d\pi_z$ and $\int g \,d\rho_z$ (resp. $\int g
\,d\pitilde_z$ and $\int g \,d\tilde{\rho}_z$) are analytic over
$\Cplus$ and are equal in $\Cnabla$ for $\im(z)$ large enough, they
are equal everywhere. Since this is true for all $g \in
C({\mathcal H})$, $\pi_z$ and $\rho_z$ (resp. $\pitilde_z$
and $\tilde{\rho}_z$) are identical on $\Cplus$.
This proves the unicity. 
\subsection{Proof of the existence of solutions}
Let us now prove the existence of solutions to (\ref{system}). Define
by recursion
$$
\pi^0(z,du,d\lambda)=\pitilde^0(z,du,d\lambda)=-\frac{1}{z} H(du,d\lambda),
$$
and
$$
\left\{
\begin{array}{ll}
\int g(u,\lambda) \pi^p(z,du,d\lambda) &= \int \frac{g(u,\lambda)}{D(\pitilde^{p-1}_z,\pi^{p-1}_z)} H(du,d\lambda)\\
\int g(u,\lambda) \pitilde^p(z,du,d\lambda) &= c \int \frac{g(cu,\lambda)}{\tilde{D}(\pi^{p-1}_z,\pitilde^{p-1}_z)}H(du,d\lambda)
+(1-c) \int_c^1 \frac{g(u,0)}{\kappa(\pi^{p-1}_z)}du
\end{array}
\right.
$$
for all $g\in C({\mathcal H})$. It is straightforward to check that 
$\pi^0_z$ (resp. $\pitilde^0_z$) is a Stieltjes kernel. 
Moreover, this remains true for $\pi^p$ and $\pitilde^p$ by induction over $p$. 
As for the unicity, we can establish
\begin{equation}\label{system-p}
\left\{
\begin{array}{ll}
\| \pi^p_z -\pi^{p-1}_z \|_{\tv} & \le {\boldsymbol \alpha}\|\pi^{p-1}_z -\pi^{p-2}_z\|_{\tv} 
+ {\boldsymbol \beta} \|\pitilde^{p-1}_z -\pitilde^{p-2}_z\|_{\tv}\\
\|\pitilde^p_z -\pitilde^{p-1}_z\|_{\tv} & \le \boldsymbol{\tilde{\alpha}}\|\pi^{p-1}_z -\pi^{p-2}_z\|_{\tv} 
+ \boldsymbol{\tilde{\beta}} \|\pitilde^{p-1}_z -\pitilde^{p-2}_z\|_{\tv}
\end{array}
\right.
\end{equation}
where $\boldsymbol{\alpha}, \boldsymbol{\tilde{\alpha}}, \boldsymbol{\beta}$ and $\boldsymbol{\tilde{\beta}}$ depend on 
$(\pi^{p-1}_z,\pitilde^{p-1}_z,\pi^{p-2}_z, \pitilde^{p-2}_z)$.  As in (\ref{majoration}), one can prove that
$$
\max(\boldsymbol{\alpha},\boldsymbol{\tilde{\alpha}}, \boldsymbol{\beta},\boldsymbol{\tilde{\beta}})\le \theta <\frac 12 
\quad \textrm{for}\ z\in \Cnabla \ \textrm{and} \ \im(z) \ \textrm{large enough}
$$ 
where $\theta$ does not depend on $(\pi^{p-1}_z,\pitilde^{p-1}_z,\pi^{p-2}_z, \pitilde^{p-2}_z)$. Therefore,
$(\pi^p_z)$ and $(\pitilde^p_z)$ are Cauchy sequences with respect to the norm $\|\cdot \|_\tv$ 
whenever $z\in \Cnabla$ and $\im(z)$ is large enough. This yields the existence and unicity 
of kernels $\pi_z$ and $\pitilde_z$
such that
$$
\forall g\in C({\mathcal H}),\quad  \left\{
\begin{array}{l}
\int g \,d\pi_z^p \xrightarrow[p\rightarrow \infty]{} \int g \,d\pi_z\\
\int g \,d\pitilde_z^p \xrightarrow[p\rightarrow \infty]{} \int g \,d\pitilde_z\\
\end{array}\right. 
$$
whenever $z\in \Cnabla$ and $\im(z)$ is large enough. 

Let $g\in C({\mathcal H})$ be fixed. Recall that $ z \mapsto \int g\,d\pi_z^p$ is analytic on $\mathbb{C}^{+}$ and that $\forall
z\in \Cplus$, $|\int g\,d\pi_z^p|\le \frac{\|g\|_{\infty}}{\im(z)}$. Therefore, $(\int g\,d\pi_z^p)_{p \geq 0}$ is a 
normal family. From every subsequence of $(\int g\,d\pi_z^p)_p$ one can thus extract a converging subsequence $(\int g\,d\pi_z^M)_M$ where
$M=M(p)$ such that 
$$
\forall K \subset \Cplus, \ K\ \mathrm{compact}\quad \sup_{z\in K}
\left| \int g\,d\pi_z^M -\Gamma(g)(z)\right| \xrightarrow[M\rightarrow \infty]{} 0,
$$
where $\Gamma(g)(z)$ is analytic over $\Cplus$. If $z\in \Cnabla$ and $\im(z)$ is large enough, 
we know that $\int g\,d\pi_z^p \rightarrow \int g\,d\pi_z$. Therefore, $\Gamma(g)(z)=\int g\,d\pi_z$ for 
$z\in \Cnabla$ and $\im(z)$ large enough. 

From this we can conclude that for all $z\in \Cplus$ , every subsequence has the same limit, say 
$\Gamma(g)(z)$ thus:
\begin{equation}\label{convpi}
\forall z\in \Cplus,\quad  \int g\,d\pi_z^p \xrightarrow[p\rightarrow \infty]{} \Gamma(g)(z),
\end{equation}
where $\Gamma(f)(z)$ is analytic. Moreover, it is straightforward to prove that
\begin{enumerate}
\item\label{linearity} $\Gamma( a\,g +b\,h)(z)= a \Gamma(g) + b \Gamma(h)$,
\item\label{continuity} $|\Gamma(g)(z)|\le \frac{\|g\|_{\infty}}{\im(z)}$,
\item $\im(\Gamma(g)(z))\ge 0$ if $g\ge 0$ and $z\in \Cplus$,
\item $\im(z \Gamma(g)(z))\ge 0$ if $g\ge 0$ and $z\in \Cplus$.
\end{enumerate}

As ${\mathcal H}$ is compact and since the application $g\mapsto
\Gamma(g)(z)$ defined for $g\in C({\mathcal H})$ is linear (property \ref{linearity}) and continuous
(property \ref{continuity}), the Riesz representation theorem yields
the existence of a measure $\pi_z(du,d\lambda)=\pi(z,du,d\lambda)$
such that
$$
\Gamma(g)(z)=\int g(u,\lambda) \pi(z,du,d\lambda).
$$
Similarly one can prove that 
\begin{equation}\label{convpitilde}
\int g\,d\ti{\pi}_z^p \rightarrow \ti{\Gamma}(g)(z)=\int g \,d\ti{\pi}_z.
\end{equation}
Let us now prove that $\pi$ and $\ti{\pi}$ satisfy
(\ref{equation1}) and (\ref{equation2}). We first check that\footnote{Il existe peut-\^etre un argument plus naturel
bas\'e sur $|d(\pi^p_z)|\ge \frac{\sqrt{2}}{2}$. A voir...}
\begin{equation}\label{denominateurnonul}
\forall z\in \Cplus, \forall u\in[0,1],\quad d(\pi_z)(u)\neq 0.
\end{equation} 

Indeed assume that for a given $u$, there exists $z_0\in \Cplus$ such that
$d(\pi_{z_0})(u)=0$ and 
consider the function $\Phi(z) = \int \sigma^2(t,cu)\,\pi(z,dt,d\zeta)$. 
As $d(\pi_z)(u)=1 + \Phi(z)$, we have $\im(\Phi(z_0)) = 0$. 
Since $\im(\Phi(z))$ is harmonic and non-negative over $\Cplus$,  the 
mean value property implies that $\im(\Phi(z)) = 0$ over $\Cplus$. 
By the Cauchy-Riemann equations, $\mathrm{Re}(\Phi(z))$ is therefore constant. But since 
$|\Phi(z)| \le \frac{\sigma_{\max}^2}{\im(z)}$, we have 
$\Phi(z) = 0$. This yields in particular $d(\pi_{z_0})(u)=1$ which contradicts 
$d(\pi_{z_0})(u)=0$.



Due to (\ref{convpi}), (\ref{convpitilde}) and
(\ref{denominateurnonul}), one has $\frac{1}{D(\ti{\pi}_z^p,\pi_z^p)(u,\lambda)} \rightarrow_p
\frac{1}{D(\ti{\pi}_z,\pi_z)(u,\lambda)}$. Since
$\frac{1}{D(\ti{\pi}_z^p,\pi_z^p)(u,\lambda)} \le \frac{1}{\im(z)}$, the
dominated convergence theorem yields:
$$
\int \frac{g(u,\lambda)}{D(\ti{\pi}_z^p,\pi_z^p)(u,\lambda)} H(du,d\lambda) \xrightarrow[p\rightarrow\infty]{} 
\int \frac{g(u,\lambda)}{D(\ti{\pi}_z,\pi_z)(u,\lambda)} H(du,d\lambda). 
$$
On the other hand $\int g\,d\pi^p_z \rightarrow_p \int g\,d\pi_z$ and (\ref{equation1}) is established.
One can establish Eq. (\ref{equation2}) similarly. 

It remains to prove that $f(z)=\int \,d\pi_z$ is the Stieltjes transform of a probability measure
(one will prove similarly the corresponding result for  $\ti f(z)=\int \,d\pitilde_z$). 
Recall that 
\begin{eqnarray*}
\im(f(z))= \im\left(\int \,d\pi_z\right)
\ge \int \frac{\im(z)}{\left|D(\pitilde_z,\pi_z)(u,\lambda)\right|^2}H(du,d\lambda)  >0
\end{eqnarray*}
by (\ref{equation1}). Moreover, since $|f(z)| \le \frac{1}{\im(z)}$, $f(z)$ is the Stieltjes transform of a 
subprobability measure. It remains to check that $\lim_{y\rightarrow +\infty} iy f(iy)=-1$.
Since 
$$
\left| \int \sigma^2(u,t)\pitilde(iy,dt,d\zeta)\right| \le \frac{\sigma_{\max}^2}{y}\quad \textrm{and}\quad  
\left| \int \sigma^2(t,cu)\,\pi(iy,dt,d\zeta)\right| \le \frac{\sigma_{\max}^2}{y},
$$
and 
$$
iy f(iy)=  \int \frac{iy\, H(du,d\lambda)}{-iy(1+\int \sigma^2(u,t)\pitilde(iy,dt,d\zeta)) +
\frac {\lambda}{1+c \int \sigma^2(t,c u) \pi(iy,dt,d\zeta)}}\ ,
$$
the Dominated convergence theorem yields the desired result.
Theorem \ref{existence-unicite} is proved.

\section{Proof of Theorem \ref{convergence}}\label{prooftheo2}

We first give an outline of the proof. 
The proof is carried out following three steps:

\begin{enumerate}
\item We first prove that for each subsequence $M(n)$ of $n$
there exists a subsequence $\Msub=\Msub(n)$ such that for all $z\in \Cplus$,
\begin{equation}\label{convergencepartout}
L^{\Msub}_{z} \xrightarrow[n\rightarrow \infty]{w} 
\mu_{z}
\quad \mathrm{and}\quad \ti L^{\Msub}_{z} \xrightarrow[n\rightarrow \infty]{w} \ti\mu_{z},
\end{equation}
where $\mu_{z}$ and $\ti\mu_{z}$ are complex measures, a priori random, with support included in ${\mathcal H}$ 
(Section \ref{step-one}).\\
\item We then prove that $z\mapsto \mu_z$ and $z\mapsto \ti\mu_z$ are Stieltjes kernels (Section \ref{step-two}).\\
\item We finally prove that for a countable collection of $z\in \Cplus$ with a limit point, say 
$$
{\mathcal C}=\{z_p\}_{p\in \mathbb{N}}\cup\{z_\infty\}\quad \textrm{with}\quad z_p\rightarrow z_\infty,
$$ 
the measures $\mu_z$ and $\ti\mu_z$ (which are a priori random) satisfy equations (\ref{equation1}) and (\ref{equation2})
almost surely for all $z\in {\mathcal C}$. Since $\, {\mathcal C}$ has a limit point in $\Cplus$, analyticity arguments will yield:
$$
\mathrm{a.s.}\quad \forall z\in \Cplus, \quad \mu_z=\pi_z \quad \mathrm{and}\quad \ti\mu_z=\pitilde_z.
$$
Otherwise stated,
$$
\mathrm{a.s.}\quad \forall z\in \Cplus,\quad L^n_z\xrightarrow[]{w}\pi_z 
\quad \mathrm{and}\quad \ti L^n_z\xrightarrow[]{w}\pitilde_z 
$$
which yields the desired result (Section \ref{step-three}).
\end{enumerate}

\subsection{Step 1: convergence of subsequences $L^{\Msub}_z$ and $\ti L^{\Msub}_z$.}\label{step-one}
Let $z_0\in \Cplus$ and let $B=\{z\in \C, |z-z_0|<\delta\} \subset \Cplus$.
Due to assumption (A-\ref{hypo3}) and to the fact that $|q_{ii}(z) |\le \im^{-1}(z)$, 
Helly's theorem implies that for each subsequence of $n$
there exists a subsequence $M=M(n)$ and a complex measure $\mu_{z_0}$ such that 
$$
 L^M_{z_0} \xrightarrow[n\rightarrow \infty]{w} \mu_{z_0}.
$$ 
Since $L^n$ is random, $\mu_{z_0}$ is a priori random too 
but due to (A-\ref{hypo3}), 
its support is included in ${\mathcal H}$. Let $(z_k,k\ge 1)$ be a sequence of complex numbers dense in $\Cplus$, then 
by Cantor diagonalization argument, one can extract a subsequence from $M$, say $\Msub$, such that  
$$
\forall k\in \mathbb{N}, \quad L^{\Msub}_{z_k} \xrightarrow[n\rightarrow \infty]{w} \mu_{z_k}
\quad \mathrm{and}\quad \ti L^{\Msub}_{z_k} \xrightarrow[n\rightarrow \infty]{w} \ti\mu_{z_k},
$$
where $\mu_{z_k}$ and $\ti\mu_{z_k}$ are complex measures, a priori random.  Let $g\in C({\mathcal K})$ 
and let $z\in \Cplus$. 
There exists $z_k$ such that $|z-z_k|\le \epsilon$ and
\begin{multline*}
\left|\int g \,dL_z^{\Msub(n)} -\int g \,dL_z^{\Msub(m)}\right| 
\le \overbrace{\left|\int g \,dL_z^{\Msub(n)} -\int g \,dL_{z_k}^{\Msub(n)}\right|}^{(a)} \\
 + \underbrace{\left|\int g \,dL_{z_k}^{\Msub(n)} -\int g \,dL_{z_k}^{\Msub(m)}\right|}_{(b)}
 + \underbrace{\left|\int g \,dL_{z_k}^{\Msub(m)} -\int g \,dL_z^{\Msub(m)}\right|}_{(c)}
\end{multline*}
Let $n$ and $m$ be large enough. Since $L_{z_k}^{\Msub}$ converges, $(b)$ goes to zero. Since $q_{ii}(z)$ 
is analytic and since $|q_{ii}(z)|\le \im^{-1}(z)$, there exists $K>0$, such that
$$
\forall i\ge 1,\ \forall z,z'\ \mathrm{close\ enough}, \quad |q_{ii}(z)-q_{ii}(z')|\le K|z-z'|.
$$
Thus $\max\{(a),(c)\} \le K \|g\|_{\infty} |z-z_k|$. Therefore,
$\left(\int g \,dL_z^{\Msub}\right)$ is a Cauchy sequence and
converges to $\Theta(g)(z)$. Since  $g\mapsto \Theta(g)(z)$ is linear and since
$|\Theta(g)(z)|\le \im^{-1}(z) \|g\|_{\infty}$, Riesz representation's theorem yields the existence of $\mu_z$ 
such that $\Theta(g)(z)=\int g\,d\mu_z$ (recall that the support of $\mu_z$ is included in ${\mathcal H}$ which is compact). 
The convergence of $\ti L_z^{\Msub}$ can be proved similarly and  (\ref{convergencepartout}) is satisfied. The first step  is proved.

\subsection{Step 2: the kernels $\mu_z$ and $\ti \mu_z$ are Stieltjes kernels.}\label{step-two} 
Let us now prove that $z\mapsto \int g\,d\mu_z$ is analytic over
$\Cplus$. Since $\left|\int g\,d L_z^{\Msub}\right|\le \im^{-1}(z)
\|g\|_{\infty}$, from each subsequence of $\left(\int g\,d
  L_z^{\Msub}\right)$, one can extract a subsequence that converges to
an analytic function. Since this limit is equal to $\int g\,d\mu_z$,
the analyticity of $z\mapsto \int g\,d\mu_z$ over $\Cplus$ is proved.
Since properties (\ref{property3}) and (\ref{property4}) defining the
Stieltjes kernels are satisfied by $L^n_z$, the kernel $\mu_z$
inherits them. Therefore, $\mu_z$ is a Stieltjes kernel. Similarly,
one can prove that $\ti\mu_z$ is a Stieltjes kernel. The second step
is proved.

\subsection{Step 3: the kernels $\mu_z$ and $\ti \mu_z$ are almost surely equal to $\pi_z$ and $\pitilde_z$}\label{step-three}
We will now prove that almost surely $\mu_z$ and $\ti\mu_z$ satisfy equations (\ref{equation1}) and (\ref{equation2}).

In the sequel we will drop the subscript $n$ from the notations relative
to matrices, and the superscript $n$ from $\Lambda_{ii}^n$. 
Let $\vec{e}_i=(\delta_{i}(k))_{1\le k\le n}$ and 
$\vec{f}_i=(\delta_{i}(k))_{1\le k\le N}$.
For the sake of simplicity, $\Sigma^T$ will be denoted $\Xi$. 
Consider the following notations:
\begin{center}
\begin{tabular}{|l|llllllll|} \hline 
Matrix &  $Y$ & $\Lambda$ & $\Sigma$ & $\Sigma_{(i)}^T$ & $Y^T$ & $\Lambda^T$ & $\Xi$ & $\Xi_{(i)}^T$\\ 
\hline \hline
$i$th row  & $\vec{y}_{i\cdot}$ & $\Lambda_{ii} \vec{e}_i$ & $\vec{\xi}_{i\cdot}$  
& $\vec{\eta}_{i\cdot}$ & $\vec{y}_{\cdot i}$& $\Lambda_{ii} \vec{f}_i$
& $\vec{\xi}_{\cdot i}$ & $\vec{\eta}_{\cdot i}$\\
\hline \hline
Matrix when $i$th row & - & - & $\Sigma_{(i)}$ & $\Sigma_{(i,i)}^T$ &-&-& $\Xi_{(i)}$ & $\Xi_{(i,i)}^T$ \\ 
is deleted &&&&&&&& \\ 
\hline 
\end{tabular}
\end{center}

In particular, 
$\vec{\xi}_{i\cdot}=\vec{y}_{i\cdot} + \Lambda_{ii} \vec{e}_i$ and 
$\vec{\xi}_{\cdot i}=\vec{y}_{\cdot i} + \Lambda_{ii} \vec{f}_i$
for $1 \leq i \leq N$. 
We will denote by $D_i$ and $\Delta_j$ the respectively 
$n\times n$ and $N\times N$ diagonal matrices defined by
$$
D_i=\textrm{diag}\left( \frac{\sigma\left(\frac iN, \frac
      1n\right)}{\sqrt{n}},\cdots,\frac{\sigma\left(\frac iN, 1\right)}{\sqrt{n}}\right)\ ,\quad
\Delta_j=\textrm{diag}\left( \frac{\sigma\left(\frac 1N, \frac
      jn\right)}{\sqrt{n}},\cdots,\frac{\sigma\left(1, \frac jn \right)}{\sqrt{n}}\right).
$$
Finally, for $1 \leq i \leq N$, we denote by 
$D_{(i,i)}$ and $\Delta_{(i,i)}$ the matrices that
remain after deleting row $i$ and column $i$ from 
$D_i$ and $\Delta_i$ respectively. 

We can state our first lemma:
\begin{lemma}\label{control} Let $z\in \Cplus$ be fixed. 
\begin{enumerate}
\item The $i$th diagonal element $q_{ii}(z)$ of the matrix 
$(\Sigma\Sigma^T -zI_N)^{-1}$ can be written~:
\begin{equation}\label{qii}
q_{ii}(z)= \frac{1}
{-z -\frac zn \sum_{k=1}^n \sigma^2\left( \frac iN, \frac kn\right) \ti{q}_{kk}(z) 
+ \frac{\Lambda_{ii}^2}
{1 +\frac 1n \sum_{k=1}^N \sigma^2\left( \frac kN, \frac in \right) q_{ii}(z) + \varepsilon_{i,n}^{(4)}+ \varepsilon_{i,n}^{(5)}}
+ \Lambda_{ii}\, \varepsilon_{i,n}^{(1)}+ \varepsilon_{i,n}^{(2)}+ \varepsilon_{i,n}^{(3)}}
\end{equation}
where  $1\le i\le N$ and 
\begin{eqnarray*}
\varepsilon_{i,n}^{(1)} &=&  -z\, \vec{y}_{i \cdot} \left( \Sigma_{(i)}^T \Sigma_{(i)} -zI\right)^{-1}\vec{e}_i^T 
- z\, \vec{e}_i \left( \Sigma_{(i)}^T \Sigma_{(i)} -zI\right)^{-1}\vec{y}_{i \cdot}^T  \\
\varepsilon_{i,n}^{(2)} &=& -z\, \vec{y}_{i \cdot} \left( \Sigma_{(i)}^T \Sigma_{(i)} -zI\right)^{-1}\vec{y}_{i \cdot}^T 
+ z\, \mathrm{Tr}\left( D_i^2\left( \Sigma_{(i)}^T \Sigma_{(i)} -zI\right)^{-1}\right) \\
\varepsilon_{i,n}^{(3)} &=& z\, \mathrm{Tr}\left( D_i^2\left( \Sigma^T \Sigma -zI\right)^{-1}\right)-
z\, \mathrm{Tr}\left( D_i^2\left( \Sigma_{(i)}^T \Sigma_{(i)} -zI\right)^{-1}\right)\\
\varepsilon_{i,n}^{(4)} &=& 
\vec{\eta}_{i \cdot} \left( \Sigma_{(i,i)} \Sigma_{(i,i)}^T -zI\right)^{-1} 
\vec{\eta}_{i \cdot}^T
- \mathrm{Tr}\left( \Delta_{(i,i)}^2
\left( \Sigma_{(i,i)} \Sigma_{(i,i)}^T -zI\right)^{-1} \right)\\
\varepsilon_{i,n}^{(5)} &=& 
\mathrm{Tr}\left( 
\Delta_{(i,i)}^2\left( \Sigma_{(i,i)} \Sigma_{(i,i)}^T -zI\right)^{-1} \right)
-\mathrm{Tr}\left( \Delta_i^2\left( \Sigma \Sigma^T -zI\right)^{-1} \right)
\end{eqnarray*}
Moreover almost surely,
\begin{equation}\label{approx}
\forall k,1\le k\le 5,\quad \lim_{n\rightarrow \infty}  \frac 1N \sum_{i=1}^N \left|\varepsilon_{i,n}^{(k)}\right| =0.
\end{equation}

\item If $1\le i\le N$ then the $i$th diagonal element 
$\ti{q}_{ii}(z)$ of the matrix $(\Sigma^T \Sigma -zI_n)^{-1}
=(\Xi \Xi^T -zI_n)^{-1}$ can be written:
\begin{equation}\label{qtilde}
\ti{q}_{ii}(z)= \frac{1}
{-z -\frac zn \sum_{k=1}^N \sigma^2\left( \frac kN, \frac in\right) q_{kk}(z) 
+ \frac{\Lambda_{ii}^2}
{1 +\frac 1n \sum_{k=1}^n \sigma^2\left( \frac iN, \frac kn \right) \ti{q}_{ii}(z) + \ti{\varepsilon}_{i,n}^{(4)}+ 
\ti\varepsilon_{i,n}^{(5)}}
+ \Lambda_{ii}\, \ti\varepsilon_{i,n}^{(1)}+ \ti\varepsilon_{i,n}^{(2)}+ \ti\varepsilon_{i,n}^{(3)}}.
\end{equation}
If $N+1\le i\le n$, then $\ti{q}_{ii}$ can be written:
\begin{equation}\label{qtildebis}
\ti{q}_{ii}(z)= \frac{1}
{-z -\frac zn \sum_{k=1}^N \sigma^2\left( \frac kN, \frac in\right) q_{kk}(z) 
+  \ti\varepsilon_{i,n}^{(2)}+ \ti\varepsilon_{i,n}^{(3)}}.
\end{equation}
where 
\begin{eqnarray*}
\ti\varepsilon_{i,n}^{(1)} &=&  -z\, \vec{y}_{\cdot i} \left( \Xi_{(i)}^T \Xi_{(i)} -zI\right)^{-1}\vec{f}_i^T 
- z\, \vec{f}_i \left( \Xi_{(i)}^T \Xi_{(i)} -zI\right)^{-1}\vec{y}_{\cdot i}^T  \\
\ti\varepsilon_{i,n}^{(2)} &=& -z\, \vec{y}_{\cdot i} \left( \Xi_{(i)}^T \Xi_{(i)} -zI\right)^{-1}\vec{y}_{\cdot i}^T 
+ z\, \mathrm{Tr}\left( \Delta_i^2\left( \Xi_{(i)}^T \Xi_{(i)} -zI\right)^{-1}\right) \\
\ti\varepsilon_{i,n}^{(3)} &=& z\, \mathrm{Tr}\left( \Delta_i^2\left( \Xi^T \Xi -zI\right)^{-1}\right)-
z\, \mathrm{Tr}\left( \Delta_i^2\left( \Xi_{(i)}^T \Xi_{(i)} -zI\right)^{-1}\right)\\
\ti\varepsilon_{i,n}^{(4)} &=& \vec{\eta}_{\cdot i} \left( \Xi_{(i,i)} \Xi_{(i,i)}^T -zI\right)^{-1} \vec{\eta}_{\cdot i}^T
- \mathrm{Tr}\left( D_{(i,i)}^2\left( \Xi_{(i,i)} \Xi_{(i,i)}^T -zI\right)^{-1} \right)\\
\ti\varepsilon_{i,n}^{(5)} &=& \mathrm{Tr}\left( D_{(i,i)}^2\left( \Xi_{(i,i)} \Xi_{(i,i)}^T -zI\right)^{-1} \right)
-\mathrm{Tr}\left( D_i^2\left( \Xi \Xi^T -zI\right)^{-1} \right)
\end{eqnarray*}
Moreover, almost surely
\begin{eqnarray}
\textrm{for}\ k&=&1,4,5\quad \lim_{n\rightarrow \infty} \frac 1N 
\sum_{i=1}^N \left|\ti\varepsilon_{i,n}^{(k)}\right| =0, \nonumber\\
\textrm{for}\ k&=&2,3\quad  \lim_{n\rightarrow \infty}\frac 1n \sum_{i=1}^n  \left|\ti\varepsilon_{i,n}^{(k)}\right| =0.
\end{eqnarray}
\end{enumerate}
\end{lemma}

\begin{proof}[Proof of Lemma \ref{control}]

Since $q_{ii}(z) = (\Sigma \Sigma^T - zI)^{-1}_{ii}$, this element is 
the inverse of the Schur complement of
$\left(\Sigma_{(i)} \Sigma_{(i)}^T -zI\right)$ in $\left(\Sigma \Sigma^T -zI\right)$ (see for instance 
\cite{Kai00}, Appendix A). In other words 
$$
q_{ii}(z) = 
\left( 
\| \vec{\xi}_{i\cdot} \|^2 - z - 
\vec{\xi}_{i\cdot} \Sigma^T_{(i)} (\Sigma_{(i)} \Sigma^T_{(i)} -zI)^{-1} 
\Sigma_{(i)} \vec{\xi}_{i\cdot}^T
\right)^{-1} \ . 
$$
Using the identity 
$$
I - \Sigma^T_{(i)} (\Sigma_{(i)} \Sigma^T_{(i)} -zI)^{-1} \Sigma_{(i)} 
=
- z (\Sigma_{(i)}^T \Sigma_{(i)} -zI)^{-1} \ , 
$$
we have 
\begin{eqnarray*}
q_{ii}(z) &=&\frac{1}{-z -z\vec{\xi}_{i\cdot} (\Sigma^T_{(i)} \Sigma_{(i)} -zI)^{-1} \vec{\xi}_{i\cdot}^T}\\
&=& \frac{1}{-z -z\vec{y}_{i\cdot} (\Sigma^T_{(i)} \Sigma_{(i)} -zI)^{-1} \vec{y}_{i\cdot}^T
-z\Lambda_{ii}^2\vec{e}_i (\Sigma^T_{(i)} \Sigma_{(i)} -zI)^{-1} \vec{e}_i^T + \Lambda_{ii}\varepsilon^{(1)}_{i,n}}\\
&=& \frac{1}{-z -\frac zn \sum_{k=1}^n \sigma^2\left(\frac iN, \frac kn\right) \ti{q}_{kk}(z) 
-z\Lambda_{ii}^2 \vec{e}_i (\Sigma^T_{(i)} \Sigma_{(i)} -zI)^{-1} \vec{e}_i^T + \Lambda_{ii}\varepsilon^{(1)}_{i,n}
+\varepsilon^{(2)}_{i,n}+\varepsilon^{(3)}_{i,n} }. 
\end{eqnarray*}
Similarly, we have 
\begin{eqnarray}\label{def-dn}
\vec{e}_i (\Sigma^T_{(i)} \Sigma_{(i)} -zI)^{-1} \vec{e}_i^T &=&  (\Sigma^T_{(i)} \Sigma_{(i)} -zI)^{-1}_{ii}\nonumber\\
&=& \frac{1}{-z -z\vec{\eta}_{i\cdot} (\Sigma_{(ii)} \Sigma_{(ii)}^T -zI)^{-1} \vec{\eta}_{i\cdot}^T}\nonumber\\
&=& \frac{1}{-z\left(1+\frac 1n \sum_{k=1}^N \sigma^2\left( \frac kN, \frac in \right) q_{kk}(z) + 
\varepsilon^{(4)}_{i,n}+\varepsilon^{(5)}_{i,n}\right) }\label{def-dn} 
\end{eqnarray}
And (\ref{qii}) is established. It is important to already note that since $\vec{\eta}_{i\cdot}$ is the $i$th row 
of $\Sigma_{(i)}^T$, $\mathbb{E} \vec{\eta}_{i\cdot}=0$ (while $\mathbb{E}\vec{y}_{i\cdot}=(0,\cdots,\Lambda_{ii},\cdots,0)$).
If $i\le N$, (\ref{qtilde}) can be established in the same way. If $i\ge N+1$, then $\vec{\xi}_{\cdot i}$ is centered:
There are no more $\Lambda_{ii}$ and all the terms involving $\Lambda_{ii}$ disappear in (\ref{qtilde}), 
which yields (\ref{qtildebis}). \\

We now prove that 
\begin{equation}\label{delta1}
 \frac 1N \sum_{i=1}^N \left|\varepsilon^{(1)}_{i,n}\right| \xrightarrow[n\rightarrow \infty]{\textrm{a.s.}}0. 
\end{equation}
One will prove similarly that $\frac 1N \sum_{i=1}^N |\ti\varepsilon^{(1)}_{i,n}|\rightarrow 0$ a.s. Denote by 
$R_n=(\Sigma^T_{(i)} \Sigma_{(i)} -zI)^{-1}=(\rho_{ij})$. Since $R_n$ is symetric, 
$\varepsilon^{(1)}_{i,n}=-z 2\vec{y}_{i\cdot} R_n \vec{e}^T_{i}$ and 
$$
|\vec{y}_{i\cdot} R_n \vec{e}^T_{i}|^4 =|\sum_{k=1}^n Y_{ik} \rho_{ki}|^4 
= \sum_{k_1,k_2,l_1,l_2} Y_{i k_1}Y_{i k_2}Y_{i l_1}Y_{i l_2} \rho_{k_1 i}\, \bar\rho_{k_2 i}\,\rho_{l_1 i}\,\bar\rho_{l_2 i}.
$$
Denote by $\mathbb{E}_{R_n}$ the expectation conditionnally to the $\sigma$-algebra generated by $R_n$.
Since $\vec{y}_{i\cdot}$ and $R_n$ are independent and since $\mathbb{E}Y_{ik}=0$, we get:
\begin{eqnarray*}
\mathbb{E}_{R_n} |\vec{y}_{i\cdot} R_n \vec{e}^T_{i}|^4 &=& 2 \mathbb{E}_{R_n} \sum_{k,l; k\neq l} Y_{ik}^2 |\rho_{ki}|^2 
Y_{il}^2 |\rho_{li}|^2\\
&\phantom{=}&+\mathbb{E}_{R_n} \sum_{k,l; k\neq l} Y_{ik}^2 \rho_{ki}^2 
Y_{il}^2 \bar\rho_{li}^2
+\mathbb{E}_{R_n} \sum_{k} Y_{ik}^4 |\rho_{ki}|^4\\
&\le & 4 \mathbb{E}(X^{n}_{ij})^{4} \frac{\sigma_{\max}^4}{n^2} \sum_{k,l} |\rho_{ki}|^2 |\rho_{li}|^2 
=  4 \mathbb{E}(X^{n}_{ij})^{4} \frac{\sigma_{\max}^4}{n^2} \left( \sum_{k} |\rho_{ki}|^2 \right)^2
\end{eqnarray*}
but 
$\sum_{k} |\rho_{ki}|^2 = 
\| R_n e_i \|^{2}  
\le \| R_n \|^2 \le \frac{1}{\im^2(z)}$.
Therefore,
\begin{equation}\label{bug-epsilon1}
\mathbb{E}\left|\varepsilon^{(1)}_{i,n}\right|^4 \le \frac{|2z|^4 4  \mathbb{E}(X^{n}_{ij})^{4} \sigma_{\max}^4}{n^2 \im^4(z)} \propto \frac{1}{n^2}.
\end{equation}
Finally,
\begin{eqnarray*}
\mathbb{P}\left\{ \frac 1N \sum_{i=1}^N \left|\varepsilon^{(1)}_{i,n}\right| > \delta \right\}
&\le & \frac{1}{\delta^4 N^4} \mathbb{E} \left( \sum_{i=1}^N |\varepsilon^{(1)}_{i,n}| \right)^4 \\
& \stackrel{(a)}{\le} & \frac{1}{\delta^4 N^4} 
\left( \sum_{i=1}^N \left( \mathbb{E}|\varepsilon^{(1)}_{i,n}|^4\right)^{\frac 14} \right)^4
\le\quad   \frac{1}{\delta^4} \sup_{1\le i\le N} \mathbb{E}|\varepsilon^{(1)}_{i,n}|^4\   \stackrel{(b)}{\propto}\quad \frac{1}{n^2} 
\end{eqnarray*}
where (a) follows from Minkowski's inequality and (b) from 
(\ref{bug-epsilon1}) and Borel-Cantelli's lemma yields Eq. (\ref{delta1}).

Let us now prove that 
\begin{equation}\label{delta2}
\frac 1N \sum_{i=1}^N \left| \varepsilon^{(2)}_{i,n}\right| \xrightarrow[n\rightarrow \infty]{\textrm{a.s.}}0. 
\end{equation}
One will prove similarly that $\frac 1n \sum |\ti\varepsilon^{(2)}_{i,n}|$, $\frac 1N \sum |\varepsilon^{(4)}_{i,n}|$ 
and $\frac 1N \sum |\ti\varepsilon^{(4)}_{i,n}|$
go to zero a.s.
Denote by $\vec{x}_{i\cdot}=(X_{i1},\cdots,X_{in})$ and write $\vec{y}_{i\cdot}=\vec{x}_{i\cdot} D_i$. In particular,
$$
\vec{y}_{i\cdot} (\Sigma^T_{(i)} \Sigma_{(i)} -zI)^{-1}\vec{y}_{i\cdot}^T =  
\vec{x}_{i\cdot} D_i (\Sigma^T_{(i)} \Sigma_{(i)} -zI)^{-1} D_i^T\vec{x}_{i\cdot}^T
$$ 
where $\vec{x}_{i\cdot}$ and $D_i (\Sigma^T_{(i)} \Sigma_{(i)} -zI)^{-1} D_i^T$ are independent. Lemma 2.7 in \cite{BaiSil98}
states that 
\begin{multline}\label{lemme-silver}
\mathbb{E}|\vec{x}_{i\cdot} C \vec{x}_{i\cdot}^T - \mathrm{Tr}\, C|^p \\\le K_p \left( 
\left( \mathbb{E}(X_{i1})^4 \mathrm{Tr}\, C C^T\right)^{p/2}
+ \mathbb{E} (X_{i1})^{2p} \mathrm{Tr}(C C^T)^{p/2} 
\right)
\end{multline}
for all $p\ge 2$. Take $p=2 + \epsilon/2$ where $\epsilon$ is given
by (A-\ref{hypo1}) and let $C=D_i (\Sigma^T_{(i)} \Sigma_{(i)}
-zI)^{-1} D_i^T$. Then 
\begin{equation}\label{major-trace}
\forall q\ge 1,\quad \mathrm{Tr}( C C^T)^q \le \frac{\sigma_{\max}^{4q}}{n^{2q-1}} \times \frac1{\im^{2q}(z)}. 
\end{equation}
Therefore, (\ref{lemme-silver}) and (\ref{major-trace}) yield 
$$
\mathbb{E}|\vec{x}_{i\cdot} C \vec{x}_{i\cdot}^T - \mathrm{Tr}\, C|^{2+\epsilon/2} 
\le \frac{K_1}{n^{1+\epsilon/4}} +\frac{K_2}{n^{1+\epsilon}} \le \frac{K}{n^{1+\epsilon/4}}
$$
where the constants $K$, $K_1$ and $K_2$ depend on the moments of
$X_{i1}$, on $\sigma_{\max}$ and on $\im(z)$. Thus
\begin{equation}\label{bug-epsilon2}
\mathbb{E}\left|\varepsilon^{(2)}_{i,n}\right|^p \le \frac{K |z|^p}{n^{1+\epsilon/4}}.
\end{equation} 
Finally,
\begin{eqnarray*}
\mathbb{P}\left\{ \frac 1N \sum_{i=1}^N \left|\varepsilon^{(2)}_{i,n}\right| > \delta \right\}
&\le & \frac{1}{\delta^p N^p} \mathbb{E} \left( \sum_{i=1}^N |\varepsilon^{(2)}_{i,n}| \right)^p \\
& \stackrel{(a)}{\le} & \frac{1}{\delta^p N^p} 
\left( \sum_{i=1}^N \left( \mathbb{E}|\varepsilon^{(2)}_{i,n}|^p\right)^{\frac 1p} \right)^p\\
&\le&   \frac{1}{\delta^p} \sup_{1\le i\le N} \mathbb{E}|\varepsilon^{(2)}_{i,n}|^p\   \stackrel{(b)}{\propto}\quad \frac{1}
{n^{1+\epsilon/4}} 
\end{eqnarray*}
where (a) follows from Minkowski's inequality and (b) 
from (\ref{bug-epsilon2}), and Borel-Cantelli's lemma yields 
(\ref{delta2}). \\ 
We now prove that \begin{equation}\label{delta3}
\frac 1N \sum_{i=1}^N \left|\varepsilon^{(3)}_{i,n}\right| \xrightarrow[n\rightarrow \infty]{\textrm{a.s.}}0. 
\end{equation}
One will prove similarly that $\frac 1n \sum |\ti\varepsilon^{(3)}_{i,n}|$ 
goes to zero.
Since $\Sigma^T \Sigma= \Sigma^T_{(i)} \Sigma_{(i)} + \vec{\xi}_{i\cdot }^T \vec{\xi}_{i\cdot }$, Lemma 2.6 in \cite{SilBai95} 
yields:
$$
\left| \mathrm{Tr} \left( (\Sigma^T \Sigma -z I)^{-1} - (\Sigma^T_{(i)} \Sigma_{(i)} -zI)^{-1}\right) D_i^2 \right|
\le \frac{\sigma_{\max}^2}{n \im(z)},
$$ 
In particular, 
\begin{equation}\label{bug-epsilon3}
\left|\varepsilon^{(3)}_{i,n}\right| \le \frac{|z| \sigma_{\max}^2}{n \im(z)} 
\end{equation}
which immediatly yields  (\ref{delta3}). \\ 
We finally prove that 
\begin{equation}\label{delta5}
\frac 1N \sum_{i=1}^N \left|\varepsilon^{(5)}_{i,n}\right| \xrightarrow[n\rightarrow \infty]{\textrm{a.s.}}0. 
\end{equation}
One will prove similarly that $\frac 1N \sum |\ti\varepsilon^{(5)}_{i,n}|$ 
goes to zero. Write
\begin{multline*}
\varepsilon_{i,n}^{(5)} = 
\mathrm{Tr}\, \Delta_{(i,i)}^2 \left( \Sigma_{(i,i)} \Sigma_{(i,i)}^T 
-zI\right)^{-1} 
-\mathrm{Tr}\, \Delta_{(i,i)}^2\left( \Sigma_{(i)} \Sigma_{(i)}^T -zI
\right)^{-1} \\
+ \mathrm{Tr}\, \Delta_{(i,i)}^2
\left( \Sigma_{(i)} \Sigma_{(i)}^T -zI\right)^{-1} 
-\mathrm{Tr}\, \Delta_i^2\left( \Sigma \Sigma^T -zI\right)^{-1} 
\end{multline*}
As for $\varepsilon_{i,n}^{(3)}$, one can prove that 
$$
\left| \mathrm{Tr}\, 
\Delta_{(i,i)}^2\left( \Sigma_{(i,i)} \Sigma_{(i,i)}^T -zI\right)^{-1} 
-\mathrm{Tr}\, 
\Delta_{(i,i)}^2\left( \Sigma_{(i)} \Sigma_{(i)}^T -zI\right)^{-1} \right| \le \frac{\sigma^2_{\max}}{n\im(z)}
$$
by applying Lemma 2.6 in \cite{SilBai95}. Let 
$$
\kappa_{i,n}=\mathrm{Tr}\, 
\Delta_{(i,i)}^2\left( \Sigma_{(i)} \Sigma_{(i)}^T -zI\right)^{-1} 
-\mathrm{Tr}\, \Delta_i^2\left( \Sigma \Sigma^T -zI\right)^{-1}.
$$
By applying to $\Sigma \Sigma^T -zI$ the identities relative to the inverse 
of a partitioned matrix (see \cite{Kai00}, Appendix A), we obtain:  
$\mathrm{Tr}\, \Delta_i^2\left( \Sigma \Sigma^T -zI\right)^{-1} 
= \Psi_1 + \Psi_2 + \Psi_3\ $ where 
\begin{eqnarray*} 
\Psi_1 &=& 
\mathrm{Tr}\, 
\Delta_{(i,i)}^2\left( \Sigma_{(i)} \Sigma_{(i)}^T -zI\right)^{-1} \\ 
\Psi_2 &=&  
\frac{\mathrm{Tr}\, 
\Delta_{(i,i)}^2
\left( \Sigma_{(i)} \Sigma_{(i)}^T -zI\right)^{-1} 
\Sigma_{(i)} \vec{\xi}_{i\cdot}^T 
\vec{\xi}_{i\cdot} \Sigma_{(i)}^T 
\left( \Sigma_{(i)} \Sigma_{(i)}^T -zI\right)^{-1} } 
{-z -z\vec{\xi}_{i\cdot} (\Sigma^T_{(i)} \Sigma_{(i)} -zI)^{-1} 
\vec{\xi}_{i\cdot}^T}  \\ 
\Psi_3 &=& 
\frac 1n 
\frac{ \sigma^2\left(\frac iN, \frac in \right) }
{-z -z\vec{\xi}_{i\cdot} (\Sigma^T_{(i)} \Sigma_{(i)} -zI)^{-1} 
\vec{\xi}_{i\cdot}^T} \ . 
\end{eqnarray*} 
In particular, 
$\kappa_{i,n}= - \Psi_2 - \Psi_3$. We have 
\begin{eqnarray*} 
| \Psi_2 |  &=&  \left| 
\frac{ 
\vec{\xi}_{i\cdot} \Sigma_{(i)}^T 
\left( \Sigma_{(i)} \Sigma_{(i)}^T -zI\right)^{-1}  
\Delta_{(i,i)}^2
\left( \Sigma_{(i)} \Sigma_{(i)}^T -zI\right)^{-1} 
\Sigma_{(i)} \vec{\xi}_{i\cdot}^T 
} 
{-z -z\vec{\xi}_{i\cdot} (\Sigma^T_{(i)} \Sigma_{(i)} -zI)^{-1} 
\vec{\xi}_{i\cdot}^T}  
\right| \\
& \leq & 
\left\| \Delta_{(i,i)} \right\|^2
\frac{ 
\left\| \left( \Sigma_{(i)} \Sigma_{(i)}^T -zI\right)^{-1} 
\Sigma_{(i)} \vec{\xi}_{i\cdot}^T \right\|^2 
} 
{\left| z + z\vec{\xi}_{i\cdot} (\Sigma^T_{(i)} \Sigma_{(i)} -zI)^{-1} 
\vec{\xi}_{i\cdot}^T \right| }  
\end{eqnarray*} 
Let
$\Sigma_{(i)} = \sum_{l=1}^{N-1} \nu_l u_l v_l^T$ be a singular value decomposition of $\Sigma_{(i)}$ where 
$\nu_l$, $u_l$, and $v_l$ are respectively the singular values, left singular
vectors, and right singular vectors of $\Sigma_{(i)}$.
Then 
$$
\left\| \left( \Sigma_{(i)} \Sigma_{(i)}^T -zI\right)^{-1} 
\Sigma_{(i)} \vec{\xi}_{i\cdot}^T \right\|^2 
= 
\sum_{l=1}^{N-1} 
\frac{\nu_l^2 \left| v_l^T \vec{\xi}_{i\cdot}^T \right|^2}
{\left| \nu_l^2 - z \right|^2} 
$$
and 
$$ 
\im\left( z + z\vec{\xi}_{i\cdot} 
(\Sigma^T_{(i)} \Sigma_{(i)} -zI)^{-1} \vec{\xi}_{i\cdot}^T \right) 
=
\im(z) \left( 
1 + 
\sum_{l=1}^{N-1} 
\frac{\nu_l^2 \left| v_l^T \vec{\xi}_{i\cdot}^T \right|^2}
{\left| \nu_l^2 - z \right|^2} 
\right) \ . 
$$
As a consequence, $\left| \Psi_2 \right| \leq \frac{\sigma_{\max}^2}{n} 
\frac{1}{\im(z)}$. Furthermore, since
$\im\left( z\vec{\xi}_{i\cdot} (\Sigma^T_{(i)} \Sigma_{(i)} -zI)^{-1} 
\vec{\xi}_{i\cdot}^T \right) \geq 0$ by Proposition \ref{proprietes-resolvante}-(3), 
we have 
$\left| \Psi_3 \right| \leq \frac{\sigma_{\max}^2}{n} 
\frac{1}{\im(z)}$. Thus,
$
\left| \varepsilon_{i,n}^{(5)} \right| \leq 
\frac{\sigma_{\max}^2}{n} \frac{3}{\im(z)}$,  
which immediatly yields (\ref{delta5}). Lemma \ref{control} is proved.
\end{proof}
Recall notation $D$ introduced at the beginning of Section 3:
$$
D(\pitilde_z,\pi_z)(u,\lambda) = -z\left(1+\int \sigma^2(u,t)\pitilde(z,dt,d\zeta)\right) +
\frac {\lambda}{1+c \int \sigma^2(t,c u) \pi(z,dt,d\zeta)}.
$$
\begin{coro}\label{coro-approx} Let $z\in \Cplus$ be fixed. Then almost surely 
\begin{equation}\label{premiere-approx}
\forall g\in C({\mathcal K}),\quad
\lim_{n\rightarrow\infty} \left| \frac{1}{N}\sum_{i=1}^{N}q_{ii}^{n} g(\Lambda_{ii}^2,i/{N})
- \frac{1}{N}\sum_{i=1}^{N}\frac{g(\Lambda_{ii}^2,i/{N})}{D(\ti \mu_z,\mu_z)\left(\frac{i}{N}\right)}\right|=0.
\end{equation}
\end{coro}

\begin{proof}[Proof of Corollary \ref{coro-approx}]
Following the notations $D$ and $d$ introduced at the beginning of Section \ref{prooftheo1}, we introduce
their empirical counterparts:
\begin{eqnarray*}
d^n(u)& =& 
1 +\frac 1n \sum_{k=1}^N \sigma^2\left( \frac kN, \frac Nn u \right) q_{ii}(z) + \varepsilon_{i,n}^{(4)}+ \varepsilon_{i,n}^{(5)}\\
D^n(u)&=&-z -\frac zn \sum_{k=1}^n \sigma^2\left( u, \frac kn\right) \ti{q}_{kk}(z) 
+ \frac{\Lambda_{ii}^2}
{d^n(u) }\\
&\phantom{=}& + \Lambda_{ii}\, \varepsilon_{i,n}^{(1)}+ \varepsilon_{i,n}^{(2)}+ \varepsilon_{i,n}^{(3)}\\
\end{eqnarray*}
Since $q_{ii}= (D^n(i/N))^{-1}$ by (\ref{qii}) and $(\Sigma^T_{(i)} \Sigma_{(i)} -zI)^{-1}_{ii}= (-zd^n(i/N))^{-1}$ by 
(\ref{def-dn}), Proposition \ref{proprietes-resolvante}-(1) yields:
\begin{equation}\label{controles-Dndn}
\frac{1}{|D^n(i/N)|}\le \frac{1}{\im(z)}\quad \textrm{and}\quad \frac{1}{|d^n(i/N)|}\le \frac{|z|}{\im(z)}.
\end{equation}
On the other hand, since $\mu_z$ and $\ti\mu_z$ are Stieltjes kernels, we have:
\begin{equation}\label{controles-Dd}
\frac{1}{|D(\ti\mu_z,\mu_z)(i/N,\Lambda_{ii}^2)|}\le \frac{1}{\im(z)}\quad 
\textrm{and}\quad \frac{1}{|d(\mu_z)(i/N)|}\le \frac{|z|}{\im(z)}.
\end{equation}
Therefore,
\begin{eqnarray*}
\lefteqn{q_{ii}^{n} - \frac{1}{D(\ti \mu_z,\mu_z)\left(\frac{i}{N},\Lambda_{ii}^2\right)}
=\frac{-z\left( \int \sigma^2( i/N,\cdot)\,d\ti L^{n}_z -  \int \sigma^2( i/N,\cdot)\,d \ti \mu_z\right)}
{D^n\left(\frac{i}{N}\right) \times D(\ti \mu_z,\mu_z)\left(\frac{i}{N},\Lambda_{ii}^2\right)}  }\\
&+&\frac{\Lambda_{ii}\, \varepsilon_{i,n}^{(1)}+ \varepsilon_{i,n}^{(2)}+ \varepsilon_{i,n}^{(3)} }
{D^n\left(\frac{i}{N}\right) \times D(\ti \mu_z,\mu_z)\left(\frac{i}{N},\Lambda_{ii}^2\right)}\\
&+& \frac{\Lambda_{ii}^2\left( \frac{N}{n} \int \sigma^2(\cdot, i/{n})\,d L^{n}_z
- c\int \sigma^2(\cdot , i/{n})\,d\mu_z\right)}
{d(\mu_z)\left(\frac{i}{N}\right)\times d^{n}\left(\frac{i}{N}\right)
\times    D^n\left(\frac{i}{N}\right) \times D(\ti \mu_z,\mu_z)\left(\frac{i}{N},\Lambda_{ii}^2\right)}\\
&+& \frac{\Lambda_{ii}^2 \left(\varepsilon_{i,n}^{(4)}+ \varepsilon_{i,n}^{(5)}\right)}
{d(\mu_z)\left(\frac{i}{N}\right)\times d^{n}\left(\frac{i}{N}\right)\times    
D^n\left(\frac{i}{N}\right) \times D(\ti \mu_z,\mu_z)\left(\frac{i}{N},\Lambda_{ii}^2\right)}
\end{eqnarray*}

Recall that the $\Lambda_{ii}$'s are assumed to be bounded (say $|\Lambda_{ii}| \le K$).
Due to (\ref{controles-Dndn}) and (\ref{controles-Dd}), we get:
\begin{eqnarray*}
\lefteqn{\left| q_{ii}^{n} - \frac{1}{D(\ti \mu_z,\mu_z)\left(\frac{i}{N},\Lambda_{ii}^2\right)}\right|
\le \frac{|z|}{\im^2(z)}\underbrace{\left|\int \sigma^2( i/N,\cdot)\,d\ti L^{n}_z 
-  \int \sigma^2( i/N,\cdot)\,d \ti \mu_z\right|}_{I(i,n)}}\\
&+& \frac{K\left|\varepsilon_{i,n}^{(1)}\right|+ \left|\varepsilon_{i,n}^{(2)}\right|+ 
\left|\varepsilon_{i,n}^{(3)}\right|}{\im^2(z)}  \\
&+& \frac{|z|^2 K^2}{\im^4(z)}\left( \underbrace{\left|\frac{N}{n} \int \sigma^2(\cdot, i/{n})\,d L^{n}_z
- c\int \sigma^2(\cdot , i/{n})\,d\mu_z\right|}_{J(i,n)} + \left|\varepsilon_{i,n}^{(4)}\right|
+\left|\varepsilon_{i,n}^{(5)}\right| \right)  \\
\end{eqnarray*}

In order to prove $\sup_{i\le N} I(i,n)\rightarrow 0$, 
recall that $C([0,1]^2) = C([0,1])\otimes C([0,1])$. In particular, $\forall \epsilon>0,$ there
exists $k\in \mathbb{N},\ g_l\in C([0,1])$ and $h_l \in C([0,1])$ for $l\le k$ such that 
$
\sup_{x,t} \left| \sigma^2(x,t) -\sum_{l=1}^k g_l(x) h_l(t)\right| \le \epsilon. 
$ 
Therefore, 
\begin{multline*}
\left| \int \sigma^2( i/N,\cdot)\,d\ti L^{n}_z -  \int \sigma^2( i/N,\cdot)\,d \ti \mu_z \right| 
\le \sup_x \left| \int \sigma^2( x,\cdot)\,d\ti L^{n}_z -  \int \sigma^2( x,\cdot)\,d \ti \mu_z\right| 
\xrightarrow[n\rightarrow\infty]{} 0
\end{multline*}
which implies that $\sup_{i\le N} |I(i,n)|$ goes to zero. One can prove similarly that $\sup_{i\le N} J(i,n)$ 
goes to zero. Therefore, 
\begin{multline*}
\left| \frac{1}{N}\sum_{i=1}^{N}\ q_{ii}^{n} g(\Lambda_{ii}^2,i/{N})
- \frac{1}{N}\sum_{i=1}^{N}\frac{g(\Lambda_{ii}^2,i/{N})}{D(\ti \mu_z,\mu_z)\left(\frac{i}{N},\Lambda_{ii}\right)}\right|\\
\le \frac{|z| \|g\|_{\infty}}{\im^2(z)} \sup_{i\le N} I(i,n) + \frac{\|g\|_{\infty}}{\im^2(z)}
\left( K \frac{1}{N} \sum_{i=1}^{N}|\varepsilon_{i,n}^{(1)}| + \frac{1}{N} \sum_{i=1}^{N}|\varepsilon_{i,n}^{(2)}|+  
\frac{1}{N} \sum_{i=1}^{N}|\varepsilon_{i,n}^{(3)}|\right)\\
+ \frac{|z|^2 K^2 \|g\|_{\infty}}{\im^4(z)} \left( \sup_{i\le N} J(i,n) + \frac{1}{N} \sum_{i=1}^{N}|\varepsilon_{i,n}^{(4)}| 
+ \frac{1}{N} \sum_{i=1}^{N}|\varepsilon_{i,n}^{(5)}|\right).
\end{multline*}
and (\ref{premiere-approx}) is proved with the help of Lemma \ref{control}. 
\end{proof}

We now come back to the proof of the third step of Theorem \ref{convergence}.
For simplicity, we will denote by $n^*=\Msub(n)$ where $\Msub$ is defined previously, by $N^*=N(n^*)$.

A direct application of the Dominated convergence theorem yields that 
$(\lambda,u) \mapsto \frac{g(\lambda,u)}{D(\ti\mu_z,\mu_z)(u,\lambda)}$ is bounded and continuous therefore (A-\ref{hypo3})
yields
\begin{equation}\label{convergence1}
\frac{1}{N^*} \sum_{i=1}^{N^*} \frac{g(\Lambda_{ii}^2,i/{N^*}) }{D(\ti\mu_z,\mu_z)(i/{N^*},\Lambda_{ii}^2)} 
\xrightarrow[n\rightarrow\infty]{}\int \frac{g(\lambda,u) }{D(\ti\mu_z,\mu_z)(u,\lambda)}H(d\lambda,du). 
\end{equation}
Moreover, 
\begin{equation}\label{convergence2}
\frac{1}{N^*} \sum_{i=1}^{N^*} g(\Lambda_{ii}^2,i/{N^*}) q_{ii}\xrightarrow[n\rightarrow\infty]{}\int g\,d\mu_z.
\end{equation}
Consider now a countable set ${\mathcal C}$ with a limit point. Since ${\mathcal C}$ is countable, 
(\ref{premiere-approx}) holds almost surely for every $z\in {\mathcal C}$ and for every $g\in C({\mathcal K})$.
Thus (\ref{convergence1}) and (\ref{convergence2}) yield that $\mu_z$ and $\ti\mu_z$ satisfy (\ref{equation1}) (and similarly
(\ref{equation2})) almost surely for all $z\in {\mathcal C}$. \\
\indent Since $\mu_z$ and $\ti\mu_z$ are Stieltjes kernels, 
one can easily prove that $z\mapsto \int \frac{g}{D(\ti\mu_z,\mu_z)}dH$ is analytic over $\Cplus$. 
Therefore, by (\ref{equation1}), the two analytic functions $z\mapsto \int g \,d\mu_z$ and 
$z\mapsto \int \frac{g}{D(\ti\mu_z,\mu_z)}dH$ 
coincide almost surely over ${\mathcal C}$ which contains a limit point. They must be equal almost surely 
over $\Cplus$. Therefore $\mu_z$ and $\ti\mu_z$ satisfy (\ref{equation1}) (and similarly
(\ref{equation2})) almost surely for all $z\in \Cplus$.\\
\indent Since $\mu$ and $\ti\mu$ are Stieltjes kernels satisfying almost surely 
(\ref{equation1}) and (\ref{equation2}), they must be almost surely equal to the unique pair
of solutions $(\pi,\pitilde)$ by Theorem \ref{existence-unicite}. In particular, $\mu$ and $\ti\mu$ are almost surely non-random. 
Thus for every subsequence 
$M=M(n)$, 
$$
\mathrm{a.s.},\quad \forall z\in \Cplus,\quad  L_z^M \xrightarrow[n\rightarrow\infty]{w} \pi_z
\quad \mathrm{and} \quad  \ti L_z^M \xrightarrow[n\rightarrow\infty]{w} \pitilde_z.
$$
Therefore, the convergence remains true for the whole sequences $L^n_z$ and $\ti L^n_z$.
Theorem \ref{convergence} is proved.

\section{Further Results and Remarks}\label{more-results}
In this section, we present two corollaries of Theorems \ref{existence-unicite} and 
\ref{convergence}. 
We will discuss the case where $\Lambda_n=0$ and the case where the variance profile $\sigma(x,y)$
is constant. These results are already well-know (\cite{BKV96,BreSil04pre,Gir90-1,Gir01a}). 

\subsection{The Centered case} 
\begin{coro}\label{centered-case}
Assume that (A-\ref{hypo1}) and (A-\ref{hypo2}) hold. Then the empirical distribution of the eigenvalues of the matrix 
$Y_n\ Y_n^T$ converges a.s. to a non-random probability measure $\mathbb{P}$ whose Stieltjes transform $f$ is given by 
$$
f(z)=\int_{[0,1]} \pi_z(dx),
$$
where $\pi_z$ is the unique Stieltjes kernel with support included in $[0,1]$ and satisfying
\begin{equation}\label{stieltjes-centered}
\forall g\in C([0,1]),\quad \int g\,d\pi_z =
\int_0^1 \frac{g(u)}{-z +\int_0^1 \frac{\sigma^2(u,t)}{1+c \int_0^1 \sigma^2(x,t)\,\pi_z(dx)}dt}du.
\end{equation}
\end{coro}
\begin{rem} In this case, one can prove that $\pi_z$ is absolutely continuous with respect to $du$, i.e.
$\pi_z(du)=k(z,u)du$ where $z\mapsto k(z,u)$ 
is analytic and $u\mapsto k(z,u)$ is continuous. Eq. (\ref{stieltjes-centered}) becomes
\begin{equation}\label{stieltjes-centered-density}
\forall u\in [0,1],\ \forall z\in \Cplus,\quad   k(u,z)=
\frac{1}{-z +\int_0^1 \frac{\sigma^2(u,t)}{1+c \int_0^1 \sigma^2(x,t)\,k(x,z)dx}dt}.
\end{equation}
Eq. (\ref{stieltjes-centered-density}) appears (up to notational differences) in \cite{Gir90-1} and in \cite{BKV96}
in the setting of Gram matrices based on Gaussian fields.
\end{rem}

\begin{proof}
  Assumption (A-\ref{hypo3}) is satisfied with $\Lambda_{ii}^n=0$ and
  $H(du,d\lambda)=du\otimes\delta_0(\lambda)$ where $du$ denotes
    Lebesgue measure on $[0,1]$. Therefore Theorems
    \ref{existence-unicite} and \ref{convergence} yield the existence
    of kernels $\pi_z$ and $\pitilde_z$ satisfying (\ref{equation1})
    and (\ref{equation2}). It is straightforward to check that in this
    case $\pi_z$ and $\pitilde_z$ do not depend on variable $\lambda$.
    Therefore (\ref{equation1}) and (\ref{equation2}) become:
\begin{equation}\label{eq-1bis}
\int g\,d\pi_z =\int \frac{g(u)}{-z(1+\int \sigma^2(u,t)\pitilde(z,dt))} du
\end{equation}
and
\begin{eqnarray}\label{eq2-bis}
\int g\,d\pitilde_z &=&c \int_{[0,1]} \frac{g(c u)}{-z(1+c \int \sigma^2(t,c u)\pi(z, dt))} du\nonumber\\
&\phantom{=}& +(1-c) \int_{[c,1]} \frac{g(u)}{-z(1+c \int \sigma^2(t,u)\pi(z,dt))} \,du \nonumber\\
&=& \int_{[0,1]} \frac{g(u)}{-z(1+c \int \sigma^2(t,u)\pi(z, dt))} du,
\end{eqnarray}
where $g\in C([0,1])$. Replacing $\int \sigma^2(u,t)\pitilde(z,dt)$ in (\ref{eq-1bis}) by the 
expression given by (\ref{eq2-bis}), one gets the following equation satisfied by $\pi_z(du)$:
$$
\int g\,d\pi_z =\int \frac{g(u)}{-z+\int \frac{\sigma^2(u,t)}{1+c \int \sigma^2(s,t)\pi(z, ds)} dt} du
$$
\end{proof}

\subsection{The non-centered case with i.i.d. entries} 

\begin{coro}\label{conv-iid}
  Assume that (A-\ref{hypo1}) and (A-\ref{hypo2}) hold where
  $\sigma(x,y)=\sigma$ is a constant function.  Assume moreover that
  $\frac 1N \sum_{i=1}^N \delta_{\Lambda_{ii}^2}\rightarrow
  H_{\Lambda}(d\lambda)$ weakly, where $H_{\Lambda}$ has a compact
  support.  Then the empirical distribution of the eigenvalues of the
  matrix $\Sigma_n \Sigma_n^T$ converges a.s. to a non-random
  probability measure $\mathbb{P}$ whose Stieltjes transform is given
  by
\begin{equation}\label{stieltjesiid}
f(z)=\int \frac{H_{\Lambda}(\,d\lambda)}{-z(1+c\sigma^2 f(z)) +(1-c)\sigma^2 +\frac{\lambda}{1+c\sigma^2 f(z)}}. 
\end{equation}
\end{coro}
\begin{rem} Eq. (\ref{stieltjesiid}) appears in \cite{BreSil04pre} in the case where $\Sigma_n= \sigma Z_n + R_n$ where 
  $Z_n$ and $R_n$ are assumed to be independent,
  $Z_{ij}^n=\frac{X_{ij}}{\sqrt{n}}$, the $X_{ij}$ being i.i.d. and
  the empirical distribution of the eigenvalues of $R_n R_n^T$
  converging to a given probability distribution. Since $R_n$ is not assumed to be 
diagonal in \cite{BreSil04pre}, the results in \cite{BreSil04pre} do not follow from Corollary \ref{conv-iid}.
\end{rem}
\begin{proof} One can build a sequence $(i/n,\Lambda_{ii}^2)$ such that 
  $\frac 1n \sum_{i=1}^n
  \delta_{(i/n,\Lambda_{ii}^2)}\xrightarrow[]{\mathcal D} du\otimes
  H_{\Lambda}(d\lambda)$.  Therefore (A-\ref{hypo3}) is satisfied with
  $H(du,d\lambda)=du\otimes H_{\Lambda}(d\lambda)$ and Theorems
  \ref{existence-unicite} and \ref{convergence} yield the existence of
  kernels $\pi_z$ and $\pitilde_z$ satisfying (\ref{equation1}) and
  (\ref{equation2}). It is straightforward to check that in this case
  $\pi_z$ and $\pitilde_z$ do not depend on variable $u$. 
  Equation (\ref{equation1})  becomes
  $$
  \int g\,d\pi_z =\int \frac{g(u,\lambda)}{-z(1+\sigma^2\int \pitilde(z,dt,d\zeta)) +
    \frac {\lambda}{1+c \sigma^2\int \pi(z,dt,d\zeta)}} H(du,d\lambda)
  $$
Let $g(u,\lambda)=1$, then (\ref{equation1}) becomes
\begin{equation}\label{eq-special}
f(z)=\int \frac{1}{-z(1+\sigma^2\ti f(z)) +
\frac {\lambda}{1+c \sigma^2f(z)}} H_{\Lambda}(d\lambda)
\end{equation}
Denote by $f_n(z)=\frac 1N \sum_1^N q_{ii}(z)$ and by $\ti
f_n(z)=\frac 1n \sum_1^N \ti{q}_{ii}(z)=\frac 1n
\mathrm{Tr}(\Sigma_n^T \Sigma_n -zI)^{-1}$.  Since $f_n(z)=\frac 1N
\mathrm{Tr}(\Sigma_n \Sigma_n^T -zI)^{-1}$ and $\ti f_n(z)=\frac 1n
\mathrm{Tr}(\Sigma_n^T \Sigma_n -zI)^{-1}$ (recall that $N\le n$),
we have $\ti f_n(z) =\frac Nn f_n(z) +\left( 1 -\frac Nn\right) \left(
  -\frac 1z\right).$ This yields $\ti f(z) = c f(z) -\frac{1-c}{z}$.
Replacing $\ti f(z)$ in (\ref{eq-special}) by this expression, we get
(\ref{stieltjesiid}).
\end{proof}

\subsection{Statement of the results in the complex case}\label{hypo-complex}
In the complex setting, Assumptions (A-\ref{hypo1})-(A-\ref{hypo3}) must be modified in the following way:
\setcounter{assump}{0}
\begin{assump}
  The random variables $(X_{ij}^n\ ;\ 1\le i\le N,\,1\le j\le n\,,\, n\ge1)$ are complex, independent and identically
  distributed. They are centered with $\mathbb{E}|X_{ij}^n|^2=1$ and
  satisfy:
$$
\exists\ \epsilon>0,\quad \mathbb{E}|X_{ij}^n|^{4+\epsilon}<\infty.
$$
\end{assump}

\begin{assump}
  The complex function $\sigma: [0,1] \times [0,1] \rightarrow \mathbb{C}$ is
  such that $|\sigma|^2$ is continuous and therefore there exist a non-negative constant
  $\sigma_{\max}$ such that
\begin{equation}
\forall (x,y)\in [0,1]^2,\quad 0\le |\sigma(x,y)|^2 \le \sigma^2_{\max}<\infty.
\end{equation}
\end{assump}
If $\Lambda_n$ is a complex deterministic $N\times n$ matrix whose non-diagonal entries are zero, assume that:
\begin{assump}
There exists a probability measure $H(\,du,d\lambda)$ over the set $[0,1]\times \mathbb{R}$ with compact support ${\mathcal H}$ such 
that
\begin{equation}
\frac1N \sum_{i=1}^N \delta_{\left( \frac iN,\,\left|\Lambda_{ii}^n\right|^2\right)}(du,d\lambda) \xrightarrow[n\rightarrow \infty]{\mathcal D} 
H(\,du,d\lambda).
\end{equation}
\end{assump}

In Theorem \ref{existence-unicite} (Eq. (\ref{equation1}) and (\ref{equation2})), one must replace $\sigma$ by its module $|\sigma|$.
The statements of Theorem \ref{convergence} and Corollary \ref{convergence-distrib} are not modified. 

\bibliography{math}

\noindent {\sc Walid Hachem},\\
Sup\'elec (Ecole Sup\'erieure d'Electricit\'e)\\
Plateau de Moulon, 3 rue Joliot-Curie\\
91192  Gif Sur Yvette Cedex, France.\\
e-mail: walid.hachem@supelec.fr\\
\\
\noindent {\sc Philippe Loubaton},\\
IGM LabInfo, UMR 8049, Institut Gaspard Monge,\\
Universit\'e de Marne La Vall\'ee, France.\\
5, Bd Descartes, Champs sur Marne, \\
77454 Marne La Vall\'ee Cedex 2, France.\\
e-mail: loubaton@univ-mlv.fr\\
\\
\noindent {\sc Jamal Najim},\\ 
CNRS, T\'el\'ecom Paris\\ 
46, rue Barrault, 75013 Paris, France.\\
e-mail: najim@tsi.enst.fr

\end{document}